# UNIFORM MARKOV RENEWAL THEORY AND RUIN PROBABILITIES IN MARKOV RANDOM WALKS[1]

By Cheng-Der Fuh

*Academia Sinica*

Let $\{X_n, n \geq 0\}$ be a Markov chain on a general state space $\mathcal{X}$ with transition probability $P$ and stationary probability $\pi$. Suppose an additive component $S_n$ takes values in the real line $\mathbf{R}$ and is adjoined to the chain such that $\{(X_n, S_n), n \geq 0\}$ is a Markov random walk. In this paper, we prove a uniform Markov renewal theorem with an estimate on the rate of convergence. This result is applied to boundary crossing problems for $\{(X_n, S_n), n \geq 0\}$. To be more precise, for given $b \geq 0$, define the stopping time $\tau = \tau(b) = \inf\{n : S_n > b\}$. When a drift $\mu$ of the random walk $S_n$ is 0, we derive a one-term Edgeworth type asymptotic expansion for the first passage probabilities $P_\pi\{\tau < m\}$ and $P_\pi\{\tau < m, S_m < c\}$, where $m \leq \infty$, $c \leq b$ and $P_\pi$ denotes the probability under the initial distribution $\pi$. When $\mu \neq 0$, Brownian approximations for the first passage probabilities with correction terms are derived. Applications to sequential estimation and truncated tests in random coefficient models and first passage times in products of random matrices are also given.

**1. Introduction.** Let $\{X_n, n \geq 0\}$ be a Markov chain on a general state space $\mathcal{X}$ with $\sigma$-algebra $\mathcal{A}$. Suppose an additive component $S_n = \sum_{t=0}^n \xi_t$ with $S_0 = \xi_0 = 0$, taking values in the real line $\mathbf{R}$, is adjoined to the chain such that $\{(X_n, S_n), n \geq 0\}$ is a Markov chain on $\mathcal{X} \times \mathbf{R}$ with

$$\begin{aligned} P\{(X_n, S_n) \in A \times (B + s) | (X_{n-1}, S_{n-1}) = (x, s)\} \\ = P\{(X_1, S_1) \in A \times B | (X_0, S_0) = (x, 0)\} \\ = P(x, A \times B), \end{aligned} \tag{1.1}$$

Received December 2001; revised November 2002.
[1]Supported in part by NSC Grant 91-2118-M-001-016.

for all $x \in \mathcal{X}, s \in \mathbf{R}, A \in \mathcal{A}$ and $B \in \mathcal{B}$ (:= Borel $\sigma$-algebra on $\mathbf{R}$). The chain $\{(X_n, S_n), n \geq 0\}$ is called a *Markov random walk*. For an initial distribution $\nu$ on $X_0$, let $P_\nu$ denote the probability measure under the initial distribution $\nu$ on $X_0$ and let $E_\nu$ denote the corresponding expectation. If $\nu$ is degenerate at $x$, we shall simply write $P_x$ $(E_x)$ instead of $P_\nu$ $(E_\nu)$. In this paper, we shall assume that $\{X_n, n \geq 0\}$ has an invariant probability $\pi$.

For $b \geq 0$, define the stopping time

$$\tau = \tau(b) = \inf\{n : S_n > b\}, \qquad \tau_+ = \tau(0). \tag{1.2}$$

In a variety of contexts, for given $m \leq \infty$ and $c \leq b$, we need to approximate the first passage probabilities

$$P_\pi\{\tau < m\}, \tag{1.3}$$

and

$$P_\pi\{\tau < m, S_m < c\}. \tag{1.4}$$

It is known that, with some proper identifications, (1.3) is the probability that the waiting time for the $(m-1)$th customer in a single server queue exceeds $b$; it is also the probability of ruin in finite time in risk theory [cf. Asmussen (1989a, b, 2000)]. The joint probability of $\tau$ and $S_m$ in (1.4) is an important ingredient to study truncated test in random coefficient models. Note that discrete time ARCH model can be defined, with some modifications, in the framework of random coefficient models; compare Bougerol and Picard (1992).

When the increments $\xi_t$ of the random walks are independent and identically distributed (i.i.d.) random variables. Siegmund (1979, 1985) and Siegmund and Yuh (1982) developed a so-called "corrected Brownian approximation" by computing correction terms in the Brownian approximation to approximate the first passage probabilities (1.3) and (1.4). In the case of a finite state ergodic Markov chain, Asmussen (1989b) derived a first-order corrected Brownian approximation for one-barrier ruin problems in risk theory, while Fuh (1997) studied one-barrier and two-barrier boundary crossing probabilities, and derived a second-order corrected Brownian approximation in Markov random walks. Arndt (1980) studied asymptotic properties of the distribution of the supremum of a random walk on a Markov chain. Malinovskii (1986) derived asymptotic expansions in the central limit theorem of (1.4) for Harris recurrent Markov chains. For a general account on ruin probabilities, the reader is referred to Asmussen (2000) and references therein.

In this paper, we study asymptotic approximations of the first passage probabilities (1.3) and (1.4) for Markov random walks on a general state space. The limiting behavior for (1.3) and (1.4) is defined as $m \to \infty$, and $b = \zeta m^{1/2}$ and $c = \gamma m^{1/2}$ for some $\gamma \leq \zeta > 0$. When a drift $\mu$ of the random walk $S_n$ is zero, we derive one-term Edgeworth type asymptotic expansions of (1.3) and (1.4). In the case of $\mu \neq 0$, we first define the conjugate transformation of the transition probability in Markov random walks and then derive corrected Brownian approximations for the first passage probabilities (1.3) and (1.4). Motivated by the approximations of (1.3) and (1.4), we study a uniform Markov renewal theorem including a rate of convergence. There are three aspects to provide the uniform Markov renewal theorem. To begin with, the condition of uniform ergodicity with respect to a given norm considered in this paper is different from the previous one [cf. Kesten (1974), Athreya, McDonald and Ney (1978), Shurenkov (1984, 1989), Alsmeyer (1994) and Fuh and Lai (2001)] and will be applied to products of random matrices. Second, we study the Markov renewal theorem with an estimate on the exponential rate of convergence. Early work in the respect can be found in Silvestrov (1994), using coupling. When the increments $\xi_t$ of the random walks are i.i.d. random variables, rates of convergence for the renewal theorem can be found in Stone (1965a, b), Carlsson and Wainger (1982), Carlsson (1983), Kartashov (1996) and Kovalenko, Kuznetsov and Shurenkov (1996). Third, the renewal theorem is in a uniform version in the sense of varying drifts. Uniform renewal theorems for simple random walks have been studied extensively in the literature; the reader is referred to Lai (1976), Kartashov (1980), Zhang (1989) and Silvestrov (1978, 1979, 1995), and references therein.

The remainder of the paper is organized as follows. In Section 2, we formulate the problem and state our main results: a uniform Markov renewal theorem with rate of convergence; one-term Edgeworth type asymptotic expansions for the first passage probabilities (1.3) and (1.4) in the case of zero drift; and corrected Brownian approximations for the first passage probabilities (1.3) and (1.4) when $\mu \neq 0$. The proofs are given in Sections 4–6, respectively. Applications to sequential estimation and truncated tests in random coefficient models and first passage times in products of random matrices are in Section 3.

**2. Main results.** Let $\{(X_n, S_n),\ n \geq 0\}$ be a Markov random walk on $\mathcal{X} \times \mathbf{R}$. For ease of notation, write $P(x, A) = P(x, A \times \mathbf{R})$ as the transition probability kernel of $\{X_n, n \geq 0\}$. For two transition probability kernels $Q(x, A), K(x, A)$, $x \in \mathcal{X}$, $A \in \mathcal{A}$ and for all measurable functions $h(x), x \in \mathcal{X}$, define $Qh$ and $QK$ by $Qh(x) = \int Q(x, dy)h(y)$ and $QK(x, A) = \int K(x, dy)Q(y, A)$, respectively.

Let $\mathcal{N}$ be the Banach space of measurable functions $h : \mathcal{X} \to \mathbf{C}$ (:= the set of complex numbers) with norm $\|h\| < \infty$. We also introduce the Banach space $\mathcal{B}$ of transition probability kernels $Q$ such that the operator norm $\|Q\| = \sup\{\|Qg\|; \|g\| \leq 1\}$ is finite. Two prototypical norms considered in the literature are the supremum norm and the $L_p$ norm. Two other commonly used norms in applications are the *weighted variation norm* and the *bounded Lipschitz norm*, described as follows:

1. Let $w : \mathcal{X} \to [1, \infty)$ be a measurable function, define for all measurable functions $h$, a weighted variation norm $\|h\|_w = \sup_{x\in\mathcal{X}} |h(x)|/w(x)$, and set $\mathcal{N}_w = \{h : \|h\|_w < \infty\}$. The corresponding norm in $\mathcal{B}_w$ is of the form $\|Q\|_w = \sup_{x\in\mathcal{X}} \int |Q|(x, dy)w(y)/w(x)$.
2. Let $(\mathcal{X}, d)$ be a metric space. For any continuous function $h$ on $\mathcal{X}$, the Lipschitz seminorm is defined by $\|h\|_{\mathrm{L}} := \sup_{x\neq y} |h(x) - h(y)|/d(x, y)$. The supremum norm is $\|h\|_\infty = \sup_{x\in\mathcal{X}} |h(x)|$. Let $\|h\|_{\mathrm{BL}} := \|h\|_{\mathrm{L}} + \|h\|_\infty$ and $\mathcal{N}_{\mathrm{BL}} = \{h : \|h\|_{\mathrm{BL}} < \infty\}$. Here BL stands for "bounded Lipschitz" and $\mathcal{N}_{\mathrm{BL}}$ is the Banach space of all bounded continuous Lipschitz functions on $\mathcal{X}$.

Denote by $P^n(x, A) = P\{X_n \in A|X_0 = x\}$, the transition probabilities over $n$ steps. The kernel $P^n$ is an $n$-fold power of $P$. Define also the Césaro averages $P^{(n)} = \sum_{j=0}^n P^j/n$, where $P^0 = P^{(0)} = I$ and $I$ is the identity operator on $\mathcal{B}$.

Definition 1. A Markov chain $\{X_n, n \geq 0\}$ is said to be uniformly ergodic with respect to a given norm $\|\cdot\|$, if there exists a stochastic kernel $\Pi$ such that $P^{(n)} \to \Pi$ as $n \to \infty$ in the induced operator norm in $\mathcal{B}$. The Markov chain $\{X_n, n \geq 0\}$ is called $w$-uniformly ergodic in the case of weighted variation norm.

The Markov chain $\{X_n, n \geq 0\}$ is assumed to be irreducible [with respect to a maximal irreducible measure $\varphi$ on $(\mathcal{X}, \mathcal{A})$], aperiodic and uniformly ergodic with respect to a given norm $\|\cdot\|$. In this paper, we assume $\varphi$ is $\sigma$-finite, and $\{X_n, n \geq 0\}$ has an invariant probability $\pi$. It is known [cf. Theorem 13.3.5 of Meyn and Tweedie (1993)] that for an aperiodic and irreducible Markov chain, if there exists some $\lambda$-small set $\mathbf{C}$ and some $P^\infty(\mathbf{C}) > 0$ such that, as $n \to \infty$,

$$\int_{\mathbf{C}} \lambda_{\mathbf{C}}(dx)(P^n(x, \mathbf{C}) - P^\infty(\mathbf{C})) \to 0,$$

where $\lambda_{\mathbf{C}}(\cdot) = \lambda(\cdot)/\lambda(\mathbf{C})$ is normalized to a probability on $\mathbf{C}$, then the chain is positive, and there exists a $\varphi$-null set $N$ such that, for any initial distribution $\nu$ with $\nu(N) = 0$,

$$\left\| \int \nu(dx)P^n(x, \cdot) - \pi \right\|_{\mathrm{tv}} \to 0 \qquad \text{as } n \to \infty,$$

where $\|\cdot\|_{\mathrm{tv}}$ denotes the total variation norm. Theorem 1.1 of Kartashov (1996) gives that $P$ has a unique *stationary projector* $\Pi$ in the sense that $\Pi^2 = \Pi = P\Pi = \Pi P$, and $\Pi(x, A) = \pi(A)$ for all $x \in \mathcal{X}, A \in \mathcal{A}$.

The following assumptions will be used throughout this paper.

C1. There exists a measure $\Psi$ on $\mathcal{X} \times \mathbf{R}$, and measurable function $h$ on $\mathcal{X}$ such that $\int \pi(dx)h(x) > 0$, $\Psi(\mathcal{X} \times \mathbf{R}) = 1$, $\int \Psi(dx \times \mathbf{R})h(x) > 0$, and the kernel $T(x, A \times B) = P(x, A \times B) - h(x)\Psi(A \times B)$ is nonnegative for all $A \in \mathcal{A}$ and, $B \in \mathcal{B}$.

C2. For all $x \in \mathcal{X}$, $\sup_{\|h\| \le 1} \|E[h(X_1)|X_0 = x]\| < \infty$.

C3. $\sup_x E_x|\xi_1|^2 < \infty$ and, for all $x \in \mathcal{X}$, $\sup_{\|h\| \le 1} \|E[|\xi_1|^r h(X_1)|X_0 = x]\| < \infty$ for some $r \ge 1$.

C4. Let $\nu$ be an initial distribution of the Markov chain $\{X_n, n \ge 0\}$; assume that for some $r \ge 1$,

$$\|\nu\| := \sup_{\|h\| \le 1} \left| \int_{x \in \mathcal{X}} h(x) E_x |\xi_1|^r \nu(dx) \right| < \infty. \tag{2.1}$$

C5. Assume that for some $n_0 \ge 1$, $\int_{-\infty}^{\infty} \int_{x \in \mathcal{X}} |E_x\{\exp(i\theta\xi_1)\}|^{n_0} \pi(dx)\, d\theta < \infty$.

C6. There exists $\Theta \subset \mathbf{R}$ containing an interval of zero such that, for all $x \in \mathcal{X}$ and $\theta \in \Theta$, $\sup_{\|h\| \le 1} \|E[\exp(\theta\xi_1)h(X_1)|X_0 = x]\| \le C < \infty$, for some $C > 0$.

C7. There exists a $\sigma$-finite measure $M$ on $(\mathcal{X}, \mathcal{A})$ such that, for all $x \in \mathcal{X}$, the probability measure $P_x$ on $(\mathcal{X}, \mathcal{A})$ defined by $P_x(A) = P(X_1 \in A|X_0 = x)$ is absolutely continuous with respect to $M$, so that $P_x(A) = \int_A p(x, y) M(dy)$ for all $A \in \mathcal{A}$, where $p(x, \cdot) = dP_x/dM$.

Remark 1. Condition C1 is a mixing condition on the Markov chain $\{(X_n, S_n), n \ge 0\}$. It is also called a minorization condition in Ney and Nummelin (1987), where they constructed a regeneration scheme and proved large deviation theorem. An alternative condition for C1 is that there exists a measure $\Psi$ on $\mathcal{X}$, and family of measures $\{h(x, B); B \in \mathcal{B}\}$ on $\mathbf{R}$, for each $x \in \mathcal{X}$ such that the kernel $T(x, A \times B) = P(x, A \times B) - h(x, B)\Psi(A)$ is nonnegative for all $A \in \mathcal{A}$ and $B \in \mathcal{B}$. If a Markov chain is Harris recurrent, then C1 holds for $n$-step transition probability. It is known that under the irreducible assumption, C1 implies that $\{(X_n, \xi_n), n \ge 0\}$ is Harris recurrent [cf. Theorem 3.7 and Proposition 3.12 of Nummelin (1984) and Theorem 4.1(iv) of Ney and Nummelin (1987)]. An example on page 9 of Kartashov (1996) also shows that there exists a uniformly ergodic Markov chain with respect to a given norm, which is not Harris recurrent. Theorem 2.2 of Kartashov (1996) states that under condition C1, a Markov chain $\{X_n, n \ge 0\}$ with transition kernel $P$ is uniformly ergodic with respect to a given norm if and only if there exists $0 < \rho < 1$ such that

$$\|P^n - \Pi\| = O(\rho^n), \tag{2.2}$$

as $n \to \infty$. When the Markov chain is uniformly ergodic with respect to the weighted variation norm, (2.2) still hold without condition C1.

Remark 2. Conditions C2–C4 are standard moment conditions. Condition C6 implies that the exponential moment, in the sense of the corresponding norm, of $\xi_1$ exists for $\theta$ in $\Theta$. Condition C5 implies that for all $n \geq n_0$, $S_n$ has a bounded probability density function for given $X_n$. The existence of the transition density in C7 will be used in Theorems 2 and 3 only. It holds in most applications.

Next, we will describe the uniform version of conditions C1–C7. Since the uniform version is in the sense of varying drift [see (2.13)], we consider a compact set $\Gamma \subset \mathbf{R}$ which contains an interval of 0. For each $\alpha \in \Gamma$, let $\{(X_n^\alpha, S_n^\alpha), n \geq 0\}$ be the Markov random walk on a general state space $\mathcal{X}$ defined as (1.1), with transition probability $P^\alpha$ and invariant probability measure $\pi^\alpha$. For each $\alpha \in \Gamma$, the Markov chain $\{X_n^\alpha, n \geq 0\}$ is assumed to be irreducible [with respect to a maximal irreducible measure $\varphi$ on $(\mathcal{X}, \mathcal{A})$], aperiodic and uniformly ergodic with respect to a given norm $\|\cdot\|$.

To establish the uniform Markov renewal theorem, we shall make use of the uniform version of (1.1) in conjunction with the following extension of the uniform Cramér's (strong nonlattice) condition:

$$g(\theta) := \inf_{\alpha \in \Gamma} \inf_{|v| > \theta} |1 - E_\pi^\alpha\{\exp(ivS_1^\alpha)\}| > 0 \qquad \text{for all } \theta > 0. \tag{2.3}$$

Additive component $S_n$ is called strongly nonlattice if $\Gamma$ has only one element. In addition, we also assume that the [conditional uniform] Cramér's (strong nonlattice) condition holds. There exists $m \geq 1$ such that

$$\sup_{\alpha \in \Gamma} \limsup_{|\theta| \to \infty} |E^\alpha\{\exp(i\theta S_m^\alpha) | X_0, X_m\}| < 1. \tag{2.4}$$

Note that under condition C5, (2.3) and (2.4) can be removed. Next, we assume the strong mixing condition holds. There exist $\gamma_1 > 0$ and $0 < \rho_1 < 1$ such that, for all $k \geq 0$ and $n \geq 1$, and for all real-valued measurable functions $g, h$ with $g, h \in \mathcal{N}$,

$$|E_\nu\{g(X_k)h(X_{k+n})\} - \{E_\nu g(X_k)\}\{E_\nu h(X_{k+n})\}| \leq \gamma_1 \rho_1^{n-1}.$$

Note that when the norm is the weighted variation norm, we only need that (2.3) and (2.4) hold without the strong mixing condition.

For $\alpha \in \Gamma$, the uniform versions of C1–C6 are:

K1. There exists a measure $\Psi^\alpha$ on $\mathcal{X} \times \mathbf{R}$, and measurable function $h$ on $\mathcal{X}$ such that $\int \pi^\alpha(dx)h(x) > 0$, $\Psi^\alpha(\mathcal{X} \times \mathbf{R}) = 1$, $\int \Psi^\alpha(dx \times \mathbf{R})h(x) > 0$, and the kernel $T^\alpha(x, A \times B) = P^\alpha(x, A \times B) - h(x)\Psi^\alpha(A \times B)$ is nonnegative for all $A \in \mathcal{A}$ and $B \in \mathcal{B}$.

K2. For all $x \in \mathcal{X}$, $\sup_{\alpha \in \Gamma} \sup_{\|h\| \leq 1} \|E^\alpha[h(X_1^\alpha) | X_0^\alpha = x]\| < \infty$.

K3. $\sup_{\alpha \in \Gamma} \sup_x E_x^\alpha |\xi_1^\alpha|^2 < \infty$ and for all $x \in \mathcal{X}$, $\sup_{\alpha \in \Gamma} \sup_{\|h\| \leq 1} \|E^\alpha[|\xi_1^\alpha|^r \times h(X_1^\alpha) | X_0^\alpha = x]\| < \infty$ for some $r \geq 1$.

K4. Let $\nu^\alpha$ be an initial distribution of the Markov chain $\{X_n^\alpha, n \geq 0\}$; assume that for some $r \geq 1$,

$$\sup_{\alpha \in \Gamma} \sup_{\|h\| \leq 1} \left| \int_{x \in \mathcal{X}} h(x) E_x^\alpha |\xi_1^\alpha|^r \nu^\alpha(dx) \right| < \infty.$$

K5. Assume that for some $n_0 \geq 1$,

$$\sup_{\alpha \in \Gamma} \int_{-\infty}^{\infty} \int_{x \in \mathcal{X}} |E_x^\alpha \{\exp(i\theta\xi_1^\alpha)\}|^{n_0} \pi^\alpha(dx)\, d\theta < \infty.$$

K6. There exists $\Theta \subset \mathbf{R}$ containing an interval of zero such that, for all $x \in \mathcal{X}$ and $\theta \in \Theta$, $\sup_{\alpha \in \Gamma} \sup_{\|h\| \leq 1} \|E^\alpha[\exp(\theta\xi_1^\alpha) h(X_1^\alpha) | X_0^\alpha = x]\| \leq C$ for some $C > 0$.

THEOREM 1. *Let $\{(X_n^\alpha, S_n^\alpha), n \geq 0\}$ be a uniformly strong nonlattice Markov random walk satisfying* K1–K4 *with $r \geq 2$ in* K3. *Let $\mu^\alpha := \mu_1^\alpha := E_\pi^\alpha \xi_1^\alpha > 0$ and $\mu_2^\alpha := E_\pi^\alpha (\xi_1^\alpha)^2 < \infty$. Then, as $s \to \infty$,*

$$\sum_{n=0}^{\infty} P_\nu^\alpha \{s \leq S_n^\alpha \leq s + h, X_n^\alpha \in A\} = \frac{h}{\mu^\alpha} \pi^\alpha(A) + o(s^{-(r-1)}), \tag{2.5}$$

$$\sum_{n=0}^{\infty} P_\nu^\alpha \{-\infty \leq S_n^\alpha \leq s, X_n^\alpha \in A\} = \left( \frac{s}{\mu^\alpha} + \frac{\mu_2^\alpha}{2(\mu^\alpha)^2} \right) \pi^\alpha(A) + o(s^{-(r-2)}). \tag{2.6}$$

*Furthermore, if* K6 *holds, then for some $r_1 > 0$, as $s \to \infty$,*

$$\sum_{n=0}^{\infty} P_\nu^\alpha \{-\infty \leq S_n^\alpha \leq s, X_n^\alpha \in A\} = \left( \frac{s}{\mu^\alpha} + \frac{\mu_2^\alpha}{2(\mu^\alpha)^2} \right) \pi^\alpha(A) + O(e^{-r_1 s}). \tag{2.7}$$

REMARK 3. When $\Gamma$ has only one element and the increments $\xi_n$ are i.i.d., these results are proved via Fourier transform and Schwartz's theory of distributions. These methods can be extended to the Markov case via perturbation theory of the transition probability operator. Such extensions can also be modified to yield the rate of convergence in Markov renewal theory, which generalizes the corresponding results of Stone (1965a, b) and Carlsson (1983) for simple random walks.

REMARK 4. For each fixed $\alpha$, the Markov renewal theorems in Theorems 1–4 in Fuh and Lai (2001) provide the rate of convergence as $s \to \infty$. However, in the applications to Theorem 3, we shall be letting $\alpha \to 0$ simultaneously with $s \to \infty$. Consequently, we must consider the possibility that certain unpleasant situations might occur, such as a case in which, as $\alpha \to 0$, the rate of convergence to 0 of the error terms in (2.5)–(2.7) gets slower and slower. Theorem 1 guarantees that this cannot happen; that is, that there is a certain rate of convergence which applies uniformly to all $\alpha$ in some neighborhood of 0.

Let $\nu$ be an initial distribution of the Markov chain $\{X_n, n \geq 0\}$ and let $\mu = E_\pi \xi_1$, $\sigma^2 = \lim_{n\to\infty} n^{-1} E_\nu\{(S_n - n\mu)^2\}$ and $\kappa = \lim_{n\to\infty} n^{-1} E_\nu\{(S_n - n\mu)^3\}$, which are well defined under C4 for some $r \geq 3$. It will be convenient to use the notation

$$P_\nu^{(m,s)}(A) = P_\nu\{A | S_m = s\}. \tag{2.8}$$

Let $\tau_+ = \inf\{n \geq 1 : S_n > 0\}$ be the first ascending ladder epoch of $S_n$, $\tau_n = \inf\{k \geq \tau_{n-1} : S_k > S_{\tau_{n-1}}\}$ be the $n$th ascending ladder epoch of $S_n$, for $n = 2, 3, \ldots,$ and let $\tau_- = \inf\{n \geq 1 : S_n \leq 0\}$ be the first descending ladder epoch of $S_n$. Since $\mu > 0$, $\tau_n$ are finite almost surely under the probability $P\{X_{\tau_+} \in A | X_0 = x\}$ and therefore, the associated ladder heights $S_{\tau_n}$ are well-defined positive random variables. Furthermore, $\{(X_{\tau_n}, S_{\tau_n}), n \geq 0\}$ is a Markov chain, and it is the so-called *ladder Markov random walk.* When $\mu = 0$, we can still define the ladder Markov chain via the property of uniform integrability in Theorem 5 of Fuh and Lai (1998). *It is assumed throughout this paper that* $P_x(\tau_+ < \infty) = 1$ *for all* $x \in \mathcal{X}$ *and that the ladder random walk is uniform ergodic with respect to a given norm.* The moment conditions C2–C4 and C6 for the ladder random walk are in Lemma 1 and Lemma 14, respectively. The uniformly strong nonlattice for the ladder random walk is in Lemma 13. Since C5 holds in Theorems 2 and 3, we do not need (2.4) anymore. Let $\pi_+$ denote the invariant measure of the kernel $P_+(x, A \times \mathbf{R}^d)$ which is *assumed to be irreducible and aperiodic.* The property of Harris recurrent for ladder Markov chains has been established in Alsmeyer (2000). Therefore, C1 holds for the ladder random walk. In Section 3, we will show how uniform ergodicity of the ladder chain and finiteness of moments of $S_{\tau_+}$ can be established in some interesting examples.

THEOREM 2. *Let* $\{(X_n, S_n), n \geq 0\}$ *be a Markov random walk satisfying* C1–C5 *and* C7 *with* $r = 3$ *in* C3. *Suppose* $\mu = 0$, $\sigma = 1$, *and that there exists* $\varepsilon > 0$ *such that* $\inf_x P_\pi\{\xi_1 \leq -\varepsilon | X_1 = x\} > 0$. *Let* $b = \zeta m^{1/2}$ *and* $s = \zeta_0 m^{1/2}$ *for some* $\zeta > 0$ *and* $-\infty < \zeta_0 < \zeta$. *Then, as* $m \to \infty$,

$$\begin{aligned} &P_\pi^{(m,s)}\{\tau < m\} \\ &\quad = \exp\{-2(b + \rho_+)(b + \rho_+ - s - \kappa/3)/(m + \kappa s/3)\} + o(m^{-1/2}), \end{aligned} \tag{2.9}$$

*where* $\rho_+ = E_{\pi_+} S_{\tau_+}^2 / 2E_{\pi_+} S_{\tau_+}$. *If in addition* $c = \gamma m^{1/2}$ *for some* $\gamma \leq \zeta$, *then, as* $m \to \infty$,

$$P_\pi\{\tau < m, S_m < c\} = \Phi\left(\frac{c + \kappa/3 - 2(b + \rho_+)}{(m + \kappa c/3)^{1/2}}\right) + o(m^{-1/2}), \tag{2.10}$$

*where* $\Phi$ *denotes the standard normal distribution function.*

REMARK 5. Approximations (2.9) and (2.10) are the corresponding results for Brownian motion with drift 0, $b$ replaced by $b+\rho_+$, $s(c)$ replaced by $s+\kappa/3$ $(c+\kappa/3)$ and $m$ replaced by $m+\kappa s/3$. Also, note that the constant $\rho_+$ in (2.9) and (2.10) reduces to $ES^2_{\tau_+}/2ES_{\tau_+}$ when $S_n$ is a simple random walk [cf. Siegmund (1985), pages 220 and 221]. Since $P_\pi\{\tau<m\}=P_\pi\{S_m\ge b\}+P_\pi\{\tau<m,S_m<b\}$, one-term Edgeworth expansion of $P_\pi\{S_m\ge b\}$ and (2.10) give a representation of (1.3).

To state Theorem 3, we need to define a twist transformation of the transition probability operator, and this requirement leads us to study the perturbation theory of certain linear operators on $\mathcal{N}$. For $z\in\mathbf{C}$, define linear operators $\mathbf{P}_z$, $\mathbf{P}$, $\nu_*$ and $\mathbf{Q}$ on $\mathcal{N}$ by

$$\begin{aligned}(\mathbf{P}_z h)(x)&=E[h(X_1)e^{z\xi_1}|X_0=x], \qquad &(\mathbf{P}h)(x)&=E[h(X_1)|X_0=x],\\ \nu_* h&=E_\nu\{h(X_0)\}, &\mathbf{Q}h&=\int h(y)\pi(dy).\end{aligned} \tag{2.11}$$

When the norm is supremum norm and $\xi_n=g(X_n)$, Nagaev (1957) and Jensen (1987) have shown that there exists sufficiently small $\delta>0$ such that, for $|z|\le\delta$, $\mathcal{N}=\mathcal{N}_1(z)\oplus\mathcal{N}_2(z)$ and

$$\mathbf{P}_z\mathbf{Q}_z h=\lambda(z)\mathbf{Q}_z h \qquad \text{for all } h\in\mathcal{N}, \tag{2.12}$$

where $\mathcal{N}_1(z)$ is a one-dimensional subspace of $\mathcal{N}$, $\lambda(z)$ is the eigenvalue of $\mathbf{P}_z$ with corresponding eigenspace $\mathcal{N}_1(z)$ and $\mathbf{Q}_z$ is the parallel projection of $\mathcal{N}$ onto the subspace $\mathcal{N}_1(z)$ in the direction of $\mathcal{N}_2(z)$. Extension of their argument to weighted variation norm and random $\xi_n$ satisfying some regularity assumptions is given in Fuh and Lai (2001). We extend this result to uniform ergodic Markov random walks with respect to a given norm in the Appendix. Let $h_1\in\mathcal{N}$ be the constant function $h_1\equiv 1$ and let $r(x;z)=(\mathbf{Q}_z h_1)(x)$. From (2.12), it follows that $r(\cdot;z)$ is an eigenfunction of $\mathbf{P}_z$ associated with the eigenvalue $\lambda(z)$, that is, $r(\cdot;z)$ generates the one-dimensional eigenspace $\mathcal{N}_1(z)$.

In particular, $z=\alpha\in\mathbf{R}$ such that there exists $\delta>0$ and $\alpha\in[-\delta,\delta]:=\Gamma$. Define the "twisting" transformation

$$P^\alpha(x,dy\times ds)=\frac{r(y;\alpha)}{r(x;\alpha)}e^{-\Lambda(\alpha)+\alpha s}P(x,dy\times ds) \qquad \text{where } \Lambda=\log\lambda. \tag{2.13}$$

Then $P^\alpha$ is the transition probability of a Markov random walk $\{(X_n^\alpha,S_n^\alpha), n\ge 0\}$, with invariant probability $\pi^\alpha$. The function $\Lambda(\alpha)$ is normalized so that $\Lambda(0)=\Lambda^{(1)}(0)=0$, where $^{(1)}$ denotes the first derivative. Then $P^0=P$ is the transition probability of the Markov random walk $\{(X_n,S_n),n\ge 0\}$ with invariant probability $\pi$. Here and in the sequel, we denote $P_\nu^\alpha$ as the probability measure of the Markov random walk $\{(X_n^\alpha,S_n^\alpha),n\ge 0\}$ with transition probability kernel (2.13), and having initial distribution $\nu^\alpha$. For ease of notation, we denote $\nu^\alpha := \nu$, and let $E^\alpha_\nu$ be the expectation under $P^\alpha_\nu$.

It is known that $\Lambda$ is a strictly convex and real analytic function for which $\Lambda^{(1)}(\alpha) = E^\alpha_\pi \xi^\alpha_1$. Therefore, $E^\alpha_\pi \xi^\alpha_1 <, =, \text{or} > 0 \Leftrightarrow \alpha <, =, \text{or} > 0$. For any value $\alpha \neq 0$ and $|\alpha| < \delta$, there is at most one value $\alpha'$ with $|\alpha'| < \delta$, necessarily of opposite sign, for which $\Lambda(\alpha) = \Lambda(\alpha')$. Assume such $\alpha'$ exists; we may let that $\alpha_0 = \min(\alpha, \alpha')$ and $\alpha_1 = \max(\alpha, \alpha')$ such that $\alpha_0 < 0 < \alpha_1$ and $\Lambda(\alpha_0) = \Lambda(\alpha_1)$. Denote $\Delta = \alpha_1 - \alpha_0$. We also assume, without loss of generality, that $\sigma^2 = \Lambda^{(2)}(0) = 1$, where ${}^{(2)}$ denotes the second derivative.

THEOREM 3. *Let $\{(X_n, S_n), n \geq 0\}$ be a strong nonlattice Markov random walk satisfying* C1–C5 *and* C7. *Let $\{(X^\alpha_n, S^\alpha_n), n \geq 0\}$ be the Markov random walk induced by* (2.13). *Suppose there exists $\varepsilon > 0$ such that $\inf_x P_\pi\{\xi_1 \leq -\varepsilon | X_1 = x\} > 0$. Let $b = \zeta m^{1/2}$ for some $\zeta > 0$, $c = \gamma m^{1/2}$ for some $\gamma \leq \zeta$, and also that $\sqrt{m}\Delta = \delta$ is a fixed positive constant. Then as $m \to \infty$, for $j = 0$ or 1,*

$$\begin{aligned} &P^{\alpha_j}_{\pi^\alpha}\{\tau < m, S^{\alpha_j}_m < c\} \\ &\quad = \int_{x \in \mathcal{X}} \frac{r(X_{\tau_+}; \alpha_j)}{r(x; \alpha_j)} \frac{r(x; \alpha_{1-j})}{r(X_{\tau_+}; \alpha_{1-j})} \\ &\qquad\qquad \times \exp[-(-1)^j \Delta(b + \rho_+)] \pi^{\alpha_j}(dx) \qquad (2.14) \\ &\qquad \times \Phi\left(\frac{c + \kappa/3 - 2(b + \rho_+)}{(m + \kappa c/3)^{1/2}} + \frac{1}{2}(-1)^j \Delta(m + kc/3)^{1/2}\right) \\ &\qquad + o(m^{-1/2}). \end{aligned}$$

REMARK 6. Let $M^{(\alpha)}_n = r(X_n; \alpha) \exp\{\alpha S_n - n\Lambda(\alpha)\}$ and $\mathcal{F}_n$ be the $\sigma$-algebra generated by $\{(X_t, S_t), t \leq n\}$. Then for $|\alpha| \leq \delta$, $\{M^{(\alpha)}_n, \mathcal{F}_n, n \geq 0\}$ is a martingale under any initial distribution $\nu$ of $X_0$ [cf. Ney and Nummelin (1987) and Fuh and Lai (1998)]. Note that $r(\cdot; \alpha)$ in (2.14) reduces to 1 when $S_n$ is a simple random walk. Hence, $r(y; \alpha)/r(x; \alpha)$ can be regarded as the reflection of Markovian dependence under uniform ergodicity condition with respect to a given norm.

**3. Examples.** In this section, we give examples of strongly nonlattice Markov random walks satisfying conditions C1–C7, and such that the underlying Markov chain $\{X_n, n \geq 0\}$ is irreducible and aperiodic. Many time series and queuing models $X_n$ are irreducible, aperiodic and $w$-uniformly ergodic Markov chains, as shown in Chapters 15 and 16 of Meyn and Tweedie (1993), and conditions C2–C4 and C6 are moment conditions on the additive components attached to $X_n$ that are satisfied in typical applications. However, renewal theorems are often applied to the ladder random walk, as in Section 2. The techniques used by Meyn and Tweedie (1993) to prove the $w$-uniform ergodicity of a rich class of time series and queuing models can also be applied to show that their ladder random walks indeed satisfy conditions C1–C7, as illustrated by the following examples.

3.1. *Random coefficient models.* Let $\{X_n, n \geq 1\}$ be the Markov chain which satisfies a first-order random coefficient autoregression model

$$X_n = \beta_n X_{n-1} + \varepsilon_n, \qquad X_0 = 0, \tag{3.1}$$

where $(\beta_n)_{n\geq 1}$ is a sequence of i.i.d. random variables with $E\beta_n = \beta$ and $\mathrm{Var}(\beta_n) = \sigma^2$, where $\sigma \geq 0$ is known. $(\varepsilon_n)_{n\geq 1}$ is a sequence of i.i.d. random variables with $E\varepsilon_n = 0$ and $\mathrm{Var}(\varepsilon_n) = 1$. Further, we assume that $(\beta_n)_{n\geq 1}$ and $(\varepsilon_n)_{n\geq 1}$ are independent, and $(\beta_n, \varepsilon_n)'$ has common density function $q$ with respect to Lebesgue measure is positive everywhere.

In the case of AR(1) model for which $\beta_n$ is a constant $\beta$, Lai and Siegmund (1983) proposed a sequential estimation procedure for the unknown parameter $\beta$. Pergamenshchikov and Shiryaev (1993) generalized their results to the model (3.1). They introduced the stopping time

$$T = T_c = \inf\left\{ n \geq 1 : \sum_{k=1}^{n} \frac{X_k^2}{1+\sigma^2 X_k^2} \geq c \right\}, \tag{3.2}$$

where $c > 0$ is a fixed number, and considered a modification of the sequential least-squares estimate

$$\widehat{b}_T = \left( \sum_{k=1}^{T} \frac{X_{k-1}X_k}{1+\sigma^2 X_{k-1}^2} \right) \Big/ \left( \sum_{k=1}^{T} \frac{X_k^2}{1+\sigma^2 X_{k-1}^2} \right). \tag{3.3}$$

Their Theorems 1–3 showed that $T_c < \infty$ with probability 1 for any $c > 0$, and $\widehat{b}_T$ is asymptotically normal under some moment conditions.

In this section we investigate the limiting behavior of $T_c$ under the stability assumption $\beta^2 + \sigma^2 < 1$. Meyn and Tweedie [(1993), Theorem 16.5.1] established $w(x) = |x|^2$-uniform ergodicity of the random coefficient model (3.1) by proving that a *drift condition* is satisfied. By Lemma 15.2.9 of Meyn and Tweedie (1993), it is also $(|x|+1)$-uniform ergodic. Suppose the conditional distribution of $\xi_n = X_n^2/(1+\sigma^2 X_n^2)$ given $X_0, \ldots, X_n$ is of the form $F_{X_{n-1},X_n}$ such that

$$\limsup_{|\theta|\to\infty} \left| \int_{\mathcal{X}} \int_{-\infty}^{\infty} \int_{-\infty}^{\infty} \left\{ \int_{-\infty}^{\infty} e^{i\theta\xi}\, dF_{x,\beta x+z}(\xi) \right\} q(\beta, z)\, d\beta\, dz\, \pi(dx) \right| < 1, \tag{3.4}$$

where $\pi$ is the stationary distribution of $\{X_n\}$. Since $\xi_1$ has probability density function with respect to Lebesgue measure, (2.4) can be removed. Let $S_n = \sum_{i=0}^n \xi_i$. Then $\{(X_n, S_n), n \geq 0\}$ is strongly nonlattice, and by Theorem 6(ii) of Fuh and Lai (1998), so is the ladder random walk with transition kernel $P_+$ defined by

$$P_+(x, A \times B) = P\{X_{\tau_+} \in A, S_{\tau_+} \in B | X_0 = x\}. \tag{3.5}$$

Assume furthermore that

$$\sup_x E_x\{\xi_1(1 + |\beta_1| + |\varepsilon_1|)\} < \infty \quad \text{and} \quad \mu := E_\pi \xi_1 > 0. \tag{3.6}$$

We first note that $X_{\tau_+}$ has a positive density function with respect to Lebesgue measure $\mathcal{L}$ and is therefore $\mathcal{L}$-irreducible. Moreover, the ladder chain is clearly aperiodic; see Section 5.4.3 of Meyn and Tweedie (1993). Let $w(x) = |x| + 1$. To show that the ladder chain is $w$-uniformly ergodic, by Theorem 16.0.1 of Meyn and Tweedie (1993), it suffices to show that there exist positive constants $b, \lambda$ and a petite set $C$ such that

$$E_x w(X_T) - w(x) \leq -\lambda w(x) + b \mathbb{1}_C(x) \qquad \text{for all } x \in \mathbf{R}, \tag{3.7}$$

for $T = \tau_+$. We first show that the drift condition (3.7) in fact holds for all stopping times $T$ [with respect to the filtration generated by $(X_n, S_n)$] such that, for some $a > 0$,

$$E_x T \leq a(|x| + 1) \qquad \text{for all} \quad x \in \mathbf{R}. \tag{3.8}$$

We then show that $\tau_+$ satisfies (3.8) and therefore (3.7) indeed holds for $T = \tau_+$.

Let $B > 0$, and denote $\varepsilon_i^* = \varepsilon_i \mathbb{1}_{\{|\varepsilon_i| \leq B\}}$, $\varepsilon_i^{**} = \varepsilon_i \mathbb{1}_{\{|\varepsilon_i| > B\}}$. Note that

$$\begin{aligned} E|X_T| &= E|\beta_T \cdots \beta_1 X_0 + \beta_{T-1} \cdots \beta_1 \varepsilon_1 + \cdots + \beta_1 \varepsilon_{T-1} + \varepsilon_T| \\ &\leq |\beta||x| + (1 - |\beta|)^{-1} B + \sum_{i=1}^T |\varepsilon_i^{**}| \qquad \text{for } X_0 = x. \end{aligned} \tag{3.9}$$

Since the $\varepsilon_i^{**}$ are i.i.d. random variables, Wald's equation yields

$$E_x \left| \sum_{i=1}^T \varepsilon_i^{**} \right| = (E_x T) E|\varepsilon_1^{**}| \leq a\ (|x| + 1) E|\varepsilon_1^{**}|, \tag{3.10}$$

by (3.8). Since $E|\varepsilon_1| < \infty$ from (3.6), we can choose $B$ sufficiently large so that $aE|\varepsilon_1^{**}| + |\beta| < 1$. In view of (3.9) and (3.10), we can then choose $0 < \lambda < 1 - |\beta| - aE|\varepsilon_1^{**}|$, $b > (1 - |\beta|)^{-1} B + aE|\varepsilon_1^{**}|$ and $C = \{x : |x| \leq K\}$ with $K$ sufficiently large such that (3.7) holds. Note that $C$ is a petite set; see Section 5.2 of Meyn and Tweedie (1993).

To show that $\tau_+$ satisfies (3.8), we use Wald's equation for Markov random walks [see Lemma 1(i) in Section 4]: For any stopping time $T$ with $E_\nu T < \infty$ and $E_\nu w(X_T) < \infty$,

$$E_\nu S_T = \mu E_\nu T + E_\nu\{\Delta(X_T) - \Delta(X_0)\}, \tag{3.11}$$

where $\sup_x |\Delta(x)|/w(x) < \infty$, and $\Delta$ is defined in (4.1). Let $\xi_i^{(B)} = \xi_i \mathbb{1}_{\{\xi_i \le B\}}$, $S_n^{(B)} = \xi_1^{(B)} + \cdots + \xi_n^{(B)}$ and $\tau(B) = \inf\{n : S_n^{(B)} > 0\}$. Since $\mu$ $(= E_\pi \xi_1 = E_\pi X_1^2/(1+\sigma^2 X_1^2)) > 0$, we can choose $B$ large enough such that $\mu^{(B)}$ $(= E_\pi \xi_1^{(B)}) > 0$. Note that $S_n^{(B)} \le S_n$ and $\tau(B) \ge \tau_+$. Hence it suffices to show that $\tau(B)$ satisfies (3.8). By the monotone convergence theorem, we need only show that (3.8) holds with $T = \tau(B) \wedge m$ for every $m \ge 1$. Since $S_{\tau(B)\wedge m} \le B$, (3.11) yields

$$\begin{aligned} B &\ge \mu^{(B)} E_x T - E_x|\Delta(X_T)| - |\Delta(x)| \\ &\ge \mu^{(B)} E_x T - cE_x|X_T| - c(|x|+2), \end{aligned} \tag{3.12}$$

since $|\Delta(x)| \le cw(x) = c(|x|+1)$ for some $c > 0$ and all $x$. By (3.9) and (3.10),

$$\begin{aligned} E_x|X_T| \le |\beta||x| &+ (1-|\beta|)^{-1} B \\ &+ (E_x T) E(|\varepsilon_1| \mathbb{1}_{\{|\varepsilon_1| > B\}}). \end{aligned} \tag{3.13}$$

Choosing $B$ large enough so that $cE|\varepsilon_1|\mathbb{1}_{\{|\varepsilon_1|>B\}} < \mu^{(B)}/2$, we obtain from (3.12) and (3.13) that

$$B + c(1-|\beta|)^{-1} B + 2c + c|x|(1+|\beta|) \ge \mu^{(B)} E_x T/2,$$

proving (3.8) for $T = \tau(B) \wedge m$. Note that the ladder chain also satisfies the mixing condition C1 by Theorem 1 of Alsmeyer (2000).

From (3.9) and (3.10) with $T = \tau_+$, it follows that $\sup_x\{E_x w(X_{\tau_+})/w(x)\} < \infty$. Hence C2 holds for the ladder chain. To prove the moment condition C6 hold for the kernel (3.5), we assume the additional moment conditions

$$\sup_x E_x \exp\{\theta(\xi_1 + \beta_1 + \varepsilon_1)\} < \infty \qquad \text{for } \theta \in \Theta. \tag{3.14}$$

First note that $E_x \exp\{\theta S_{\tau_+}\} w(X_{\tau_+}) \le \{E_x \exp\{\theta p S_{\tau_+}\}\}^{1/p} \{E_x w^q(X_{\tau_+})\}^{1/q}$, where $p^{-1} + q^{-1} = 1$. Lemma 14 in Section 6 implies that $E_x \exp\{\theta p S_{\tau_+}\} < \infty$ for $x \in \mathcal{X}$. From (3.9), there exists a constant $c_q$ depending only on $q$ such that

$$E_x|X_{\tau_+}|^q \le c_q \biggl\{ |x|^q + (1-|\beta|^{-1}) E_x \max_{i \le \tau_+} |\varepsilon_i|^q \biggr\}$$

$$\leq c_q\left\{|x|^q + (1-|\beta|^{-1})E_x\sum_{i=1}^{\tau_+}|\varepsilon_i|^q\right\}$$

$$= c_q\{|x|^q + (1-|\beta|^{-1})E|\varepsilon_1|^q E_x\tau_+\}.$$

It follows that $\sup_x E_x(\exp\{\theta S_{\tau_+}\}w(X_{\tau_+}))/w(x) < \infty$. By using the same argument and applying Wald's equation for Markov chain in Lemma 1, it follows that the moment condition C3 also holds for the ladder chain. Since $S_0 = 0$ and $\sup_x E_x \exp\{\theta S_{\tau_+}\}/w(x) < \infty$, C4 also holds for the kernel $P_+$ if the initial distribution $\nu$ satisfies $\int_{-\infty}^{\infty}|x|\,d\nu(x) < \infty$. Hence C1–C6 are satisfied by the ladder random walk with transition kernel $P_+$ when the underlying chain is the random coefficient model (3.1) and $\{\xi_n\}$ satisfies (3.4), (3.6) and (3.14).

Under the normality assumption on $(\beta_k, \varepsilon_k)$ with known $\sigma^2$, the log-likelihood ratio statistic $Z_n$ for testing $H_0\colon\beta \leq \mu_0$ against $H_1\colon\beta > \mu_0$ is given by

$$Z_n = \tfrac{1}{2}\sum_{i=1}^{n}((X_i - \mu_1 X_{i-1})^2 - (X_i - \mu_0 X_{i-1})^2), \tag{3.15}$$

where $\mu_1$ is so chosen $\mu_1 > \mu_0$. Define the stopping time $T_\lambda = \inf\{n \geq 1\colon Z_n \geq \lambda\}$. Given $m > 0$, we consider the test of $H_0\colon\beta \leq \mu_0$ against $H_1\colon\beta > \mu_0$ defined by the following: stop sampling at $\min(T_\lambda, m)$; reject $H_0$ if $T_\lambda \leq m$, and otherwise do not reject $H_0$.

Under the stability assumption $\beta^2 + \sigma^2 < 1$, and consider $Y_n = (X_{n-1}, X_n)$ as the underlying Markov chain. Suppose the conditional distribution of $\xi_n = (X_n - \beta_1 X_{n-1})^2 - (X_n - \beta_0 X_{n-1})^2$ given $Y_0, \ldots, Y_n$ is of the form $F_{Y_{n-1},Y_n}$ such that

$$\limsup_{|\theta|\to 0}\left|\int_{\mathcal{X}\times\mathcal{X}}\int_{-\infty}^{\infty}\int_{-\infty}^{\infty}\left\{\int_{-\infty}^{\infty} e^{i\theta\xi}\,dF_{y,\beta y+z}(\xi)\right\} \times q(\beta,z)\,d\beta\,dz\,\pi_1(dy)\right| < 1, \tag{3.16}$$

where $\pi_1$ is the stationary distribution of $\{Y_n\}$. Let $S_n = \sum_{i=0}^n \xi_i$; then $\{(Y_n, S_n), n \geq 0\}$ is strongly nonlattice. It is easy to see that there exists $\varepsilon > 0$ such that $\inf_x P_{\pi_1}\{\xi_1 \leq -\varepsilon | X_1 = x\} > 0$. The rest of the argument is the same as (3.4)–(3.14) and is omitted.

3.2. *Products of random matrices.* In this section, we apply the renewal theory in Section 2 to generalize some results of Kesten (1973) on products of random matrices in three directions. First, while Kesten considered products of i.i.d. matrices $M_n$, we work with the more general setting in which $\{(X_n, M_n), n \geq 0\}$ are products of Markov random matrices. Second, we provide uniform renewal theory, with polynomial and exponential rate of convergence, respectively. This extension enables us to apply corrected diffusion approximation for the first passage probabilities. Third, while Kesten assumed the entries of $M_n$ to be positive with probability 1, we can dispense with this assumption. Moreover, our proof is considerably simpler and provides a more transparent description of the basic constants that appear in his results.

We shall consider $k \times k$ nonsingular matrices $M$, with real entries, and define the norm by $\|M\| = \sup_{|x|=1} |Mx|$, where $|\cdot|$ is a norm in $\mathbf{R}^k$. Following Bougerol (1988), define $\chi(M) = \max(\log \|M\|, \log \|M^{-1}\|)$. Let $\{(X_n, M_n), n \geq 0\}$ be a Markov chain satisfying the following assumptions:

(A1) $M_n$ is a $k \times k$ nonsingular matrix with real entries such that

$$\sup_x E\{\exp(a\chi(M_1))|X_0 = x\} < \infty \qquad \text{for some } a > 0.$$

(A2) $\{X_n, n \geq 0\}$ is a $w$-uniformly ergodic Markov chain and satisfies C7.

(A3) $\{(X_n, M_n), n \geq 0\}$ is quasi-irreducible and $\gamma_1 \neq \gamma_2$, where $\gamma_1 \geq \gamma_2 \geq \cdots \geq \gamma_k$ denotes its Lyapunov exponents.

For the definition of "quasi-irreducibility," see Bougerol [(1988), page 199]. For the definition and basic properties of Lyapunov exponents, see Bougerol and Lacroix [(1985), Sections III.5 and III.6) and Bougerol [(1988), pages 197 and 198]. Let $M_0$ be the $k \times k$ identity matrix and define the product

$$\Pi_n = M_n \cdots M_0. \tag{3.17}$$

Let $u$ be unit column vectors in $\mathbf{R}^k$. Since $\Pi_n$ is nonsingular, $\Pi_n u \neq 0$ and

$$\log |\Pi_n u| = \sum_{t=1}^{n} \xi_t \qquad \text{where } \xi_t = \log(|\Pi_t u|/|\Pi_{t-1} u|). \tag{3.18}$$

Let $Y_0 = (X_0, u), \ldots, Y_n = (X_n, \Pi_n u/|\Pi_n u|)$. Define $\xi_t$ as (3.18), and let $S_n = \sum_{t=1}^n \xi_t$. Then it follows from (3.17) that $\{(Y_n, S_n), n \geq 0\}$ is a Markov random walk. Under (A1)–(A3), Bougerol (1988) has shown that $Y_n$ has an invariant measure and is uniform ergodic with Hölder continuous norm [see Definitions 3 and 4 of Bougerol (1988)]. Moreover, in view of (A1), conditions C3 and C4 are satisfied for every $r > 1$. Furthermore, condition C6 holds. Assuming $\xi_1$ to be strongly nonlattice and conditional strongly nonlattice, we can therefore apply the renewal theorems in Section 2 to the Markov random walk $\{(Y_n, S_n), n \geq 0\}$, thereby both generalizing Kesten's (1973) renewal theory for products of i.i.d. matrices with positive entries and providing convergence rates in the renewal theorems. The mean value $\mu$ in these theorems is equal to $\gamma_1$ of upper Lyapunov exponents.

Let $\mathcal{S}$ denote the sphere consisting of unit column vectors in $\mathbf{R}^k$. For $u \in \mathcal{S}$, define the stopping time

$$N(b) = \inf\{n \geq 1 : |\Pi_n u| > e^b\} = \inf\{n \geq 1 : S_n > b\}, \qquad \inf \phi = \infty. \tag{3.19}$$

Suppose $\gamma_1 \geq 0$. Since $\sup_{x\in\mathcal{X},u\in\mathcal{S}} E(|\xi_1|^r|X_0 = x) < \infty$ in view of (A1), Theorem 2 in Section 2 can be applied to show that (2.9) and (2.10) hold for $N(b)$ as $b \to \infty$. This generalizes Theorem 2 of Kesten (1973) that considers the case of i.i.d. $M_n$ with positive entries.

We next consider the case $\gamma_1 < 0$ and assume in addition that, for some $0 < a^* < a$ [where $a$ is given in (A1)],

$$\text{(A4)} \qquad \inf_{x\in\mathcal{X},u\in\mathcal{S}} E\{|M_1 u|^{a^*}|X_0 = x\} \geq k^{a^*/2}.$$

For i.i.d. matrices $M_n = (M_n(h,i))_{1\leq h,i\leq k}$ with positive entries, Kesten's (1973) Theorem 3 restricts $u$ to the subset $\mathcal{S}_+$ of $\mathcal{S}$ consisting only of vectors with nonnegative entries and assumes, among other conditions, that

$$E\left\{\min_{1\leq i\leq k}\left(\sum_{h=1}^{k} M_1(h,i)\right)^{a^*}\right\} \geq k^{a^*/2}, \tag{3.20}$$

which implies (A4) with $\mathcal{S}$ replaced by $\mathcal{S}_+$ [see the inequality preceding (2.66) in Kesten (1973)]. Under assumptions (A1) and (A2), if the assumption of quasi-irreducibility in (A3) is strengthened into "strong irreducibility" [see Section 5 of Bougerol (1988) for its definition and properties], then it can be shown that tilting for the operator $\mathbf{P}_\alpha$ is defined by

$$(\mathbf{P}_\alpha f)(x,u) = E\{e^{\alpha\xi_1} f(Y_1)|Y_0 = (x,u)\}$$

on the space $L(\alpha)$ of functions on $\mathcal{X} \times \mathcal{S}$ with the Hölder continuous norm whose spectrum is taken over $x \in \mathcal{X}$ [if the induced Markov chain $\{(Y_n, S_n), n \geq 0\}$ is irreducible]. By making use of (A1)–(A4) and Proposition 1 in the Appendix, $(\mathbf{P}_\alpha f)(x,u)$ is well defined. Let $\lambda(\alpha)$ and $r((x,u);\alpha)$ be the largest eigenvalue and associated eigenfunction defined as (2.11)–(2.13).

Therefore, the usual tilting argument shows that Theorem 3 of Kesten (1973) holds. In particular, taking $B = e^b$ yields

$$BP_x\left\{\max_{n\geq 1}|\Pi_n u| > B\right\} \longrightarrow Kr((x,u),\alpha) \qquad \text{as } B \to \infty, \tag{3.21}$$

where $K = (e^{-b(\Lambda(\alpha)-1)}/\gamma_1)\int_{\mathcal{X}\times\mathbf{R}^k} 1/r((x,u);\alpha)\,d\pi_+^\alpha(x,u)$.

**4. Proof of Theorem 1.** The ingredients we need to make the *uniform* Markov renewal theorem over the family $\{(X_n^\alpha, S_n^\alpha), n \geq 0 : \alpha \in \Gamma\}$ are provided by Lemmas 3 and 4. The proof of Lemma 3 depends on a uniform upper bound for the expectation of the overshoot, which we state and prove in Lemma 2. To prove Lemma 2, we need Wald's equation for Markov chains, and moment convergence of the stopping time $\tau(b)$ defined in (1.2) for all $b \geq 0$. These are included in Lemma 1. A version of Wald's equations for uniformly ergodic Markov random walks can be found in Fuh and Lai (1998), where they applied the spectral theory of positive operators related to Markov semigroups. Fuh and Zhang (2000) first derived Poisson equations for Markov random walks, and then applied them to establish Wald's equations. Here in C1 and C2, we applied results in Harris recurrent Markov random walks to obtain (first-order) Wald's equation via Poisson equation.

LEMMA 1. *Assume* C1 *and* C2 *and that* $0 < \mu := E_\pi \xi_1 < \infty$. *Let* $\nu$ *be an initial distribution of* $X_0$, *and let* $T$ *be a stopping time such that* $E_\nu T < \infty$:

(i) *If* $\sup_x E_x(|\xi_1|) < \infty$, *then*

$$E_\nu S_T = \mu E_\nu T + E_\nu\{\Delta(X_T) - \Delta(X_0)\}.$$

*The constant* $E_\nu\{\Delta(X_T) - \Delta(X_0)\}$ *is zero when* $\nu = \pi$. *Denote* $\xi_1^+$ ($\xi_1^-$) *as the positive (negative) part of* $\xi_1$.

(ii) *If* $\sup_x E_x(\xi_1^-) < \infty$, *then* $E_\nu \tau(b) < \infty$.

(iii) *Let* $p \geq 1$. *If* $\sup_x E_x(\xi_1^-) < \infty$ *and* $\sup_x E_x(\xi_1^+)^p < \infty$, *then* $E_\nu S^p_{\tau(b)} < \infty$.

PROOF. (i) The minorization condition C1 ensures that $\{(X_n, S_n), n \geq 0\}$ is a split chain [cf. Lemma 3.1 of Ney and Nummelin (1987)]. Under the irreducible assumption, it is also a Harris recurrent Markov chain. Proposition 17.4.1 and Theorem 17.4.2 of Meyn and Tweedie (1993) give that the following Poisson equation

$$E_x \Delta(X_1) - \Delta(x) = E_x \xi_1 - E_\pi \xi_1 \tag{4.1}$$

has a solution $\Delta : \mathcal{X} \to \mathbf{R}$ for almost every $x \in \mathcal{X}$. Under the assumption of $\sup_x E_x(|\xi_1|) < \infty$, $E_\nu(\Delta(X_T) - \Delta(X_0))$ is finite. Therefore, by Corollary 1 and Theorem 4 of Fuh and Zhang (2000), we have the proof of (i).

To prove (ii), we choose $B > 0$ such that $\mu' := E_\pi(\xi_1(B)) > 0$, where $\xi_t(B) = \xi_t I(\xi_t \leq B)$. Let $S'_n = \xi_1(B) + \cdots + \xi_n(B)$ and let $N_b = \inf\{n \geq 1 : S'_n \geq b\}$. Then $S'_n \leq S_n$ and $N_b \geq \tau(b)$. For $m > 0$, apply (i) to $N_b \wedge m$; we have $E_\nu S'_{N_b \wedge m} = \mu' E_\nu(N_b \wedge m) + O(1)$ as $m \to \infty$. By the monotone convergence theorem, $\lim_{m\to\infty} E_\nu(N_b \wedge m) = E_\nu N_b$. Moreover, by the definition of $N_b$, $S_{N_b \wedge m} \leq b + B$ for all $m \geq 1$. Hence $b + B \geq \mu' E_\nu N_b - a$ for some $a > 0$, and therefore $\infty > E_\nu N_b \geq E_\nu \tau(b)$.

Finally, we prove (iii). Since $0 \leq S_{\tau(b)} < b + \xi_{\tau(b)}$, it follows from Minkowski's inequality that

$$(E_\nu S^p_{\tau(b)})^{1/p} \leq b + \left\{ E_\nu \left[ \sum_{t=1}^{\tau(b)} (\xi_t^+)^p \right] \right\}^{1/p}.$$

Since we have already shown that $E_\nu \tau(b) < \infty$ and $\sup_x E_x(\xi_1^+)^p < \infty$ by assumption, it follows from (i) that

$$E_\nu \left[ \sum_{t=1}^{\tau(b)} (\xi_t^+)^p \right] \leq \left\{ \sup_x E_x(\xi_1^+)^p \right\} E_\nu \tau(b) + O(1),$$

proving the finiteness of $E_\nu S^p_{\tau(b)}$. $\square$

LEMMA 2. *Assume* C1–C4 *with* $r = 2$ *in* C3. *Suppose* $\mu > 0$. *For* $b \geq 0$, *let* $R(b) = S_{\tau(b)} - b$. *Then,*

$$\sup_{b\geq 0} E_\pi R(b) \leq \frac{E_\pi(\xi_1^+)^2}{E_\pi \xi_1}. \tag{4.2}$$

*When the initial distribution of* $X_0$ *is* $\nu$, (4.2) *becomes that there exists a constant* $K > 0$ *such that* $\sup_{b\geq 0} E_\nu R(b) \leq E_\pi(\xi_1^+)^2/E_\pi\xi_1 + K$.

REMARK 7. In the case of simple random walks, the upper bound (4.2) was given in Lorden (1970) by pathwise integration.

PROOF OF LEMMA 2. For any values of $\xi_1, \xi_2, \ldots$, the overshoot function $\{R(b); b \geq 0\}$ is piecewise linear, with all pieces having slope $-1$. We consider first the case where the $\xi$'s are nonnegative. It is easy to see that, for $c \geq 0$,

$$\int_0^c R(b)\,db = \int_0^{S_{\tau(c)}} R(b)\,db - \int_c^{S_{\tau(c)}} R(b)\,db = \tfrac{1}{2}\sum_{t=1}^{\tau(c)} \xi_t^2 - \tfrac{1}{2}R(c)^2. \tag{4.3}$$

Since for $c \geq 0$, $E_\pi\tau(c)$ is finite by Lemma 1(ii), the sum in (4.3) has finite expectation by C3 and Wald's equation for Markov random walks, and since the other terms are nonnegative, they also have finite expectations. Since $R(b) \geq 0$ for all $b$, we have by Fubini's theorem and Wald's equation for Markov random walks in Lemma 1(i)

$$\int_0^c E_\pi R(b)\,db = \tfrac{1}{2}E_\pi\xi_1^2 E_\pi\tau(c) - \tfrac{1}{2}E_\pi R(c)^2,$$

where $E_\pi\xi_1^2$ is finite via condition C3. Note that $E_\pi\{\Delta(X_{\tau(c)}) - \Delta(X_0)\} = 0$ in the Wald's equation.

By Jensen's inequality and Wald's equation for Markov random walks,

$$\int_0^c E_\pi R(b)\,db \leq \tfrac{1}{2}\mu^{-1}E_\pi\xi_1^2(c + E_\pi R(c)) - \tfrac{1}{2}(E_\pi R(c))^2. \tag{4.4}$$

It is easy to see that for all $b, u \geq 0$, $E_\pi\tau(b+u) \leq E_\pi\tau(b) + E_\pi\tau(u)$, since the conditional expectation of $\tau(b+u) - \tau(b)$ given $\tau(b) = n$, $X_0, X_1, \xi_1, \ldots, X_n, \xi_n$ equals $E_\pi\tau(u-r)$, where $r = \xi_1 + \cdots + \xi_n - b > 0$ and $\tau(u-r)$ is zero if $r > u$, so that $E_\pi\tau(u-r) \leq E_\pi\tau(u)$. It follows from Wald's equation for Markov random walks that $E_\pi R(b)$ is a subadditive function of $b$ and therefore

$$\begin{aligned} &\tfrac{1}{2}c(E_\pi R(c) + g_c) \\ &\quad\leq \tfrac{1}{2}c \inf_{0\leq b\leq 1/2c}(E_\pi R(b) + E_\pi R(c-b)) \\ &\quad\leq \int_0^{1/2c}(E_\pi R(b) + E_\pi R(c-b))\,db = \int_0^c E_\pi R(b)\,db. \end{aligned} \tag{4.5}$$

Combining (4.4) and (4.5) and rewriting, we obtain

$$(E_\pi R(c))^2 + (c - E_\pi \xi_1^2/\mu) E_\pi R(c) - c E_\pi \xi_1^2/\mu \le 0. \tag{4.6}$$

The left-hand side of (4.6) is a quadratic in $E_\pi R(c)$ which is nonpositive only between its roots, $-c$ and $E_\pi \xi_1^2/\mu$. Therefore, $E_\pi R(c) \le E_\pi \xi_1^2/\mu$ and since $c$ is arbitrary, the proof is complete for the nonnegative case.

The case where $\xi_1, \xi_2, \ldots$ may take negative values reduces to the nonnegative case through consideration of the associated sequence of positive ladder variables, which forms the ladder Markov random walks $\{(X_{\tau_n}, S_{\tau_n}), n \ge 0\}$ defined in the paragraph before Theorem 2. We first note that, by assumption, the ladder Markov chain is uniformly ergodic with respect to a given norm. Next, we need to verify that conditions C1–C4 still hold for the associated ladder Markov chains $\{(X_{\tau_n}, S_{\tau_n}), n \ge 0\}$. It is known [cf. Theorem 1 of Alsmeyer (2000)] that if $\{(X_n, S_n), n \ge 0\}$ is Harris recurrent, then the associated ladder Markov chains $\{(X_{\tau_n}, S_{\tau_n}), n \ge 0\}$ are also Harris recurrent. Under the irreducible assumption, the minorization condition C1 is equivalent to Harris recurrent. Therefore the ladder Markov chain $\{(X_{\tau_n}, S_{\tau_n}), n \ge 0\}$ satisfies C1. The moment conditions C2–C4 hold by Lemma 1(iii).

Since $R(b)$ is pointwise the same for $\xi_1, \xi_2, \ldots$ and the sequence of ladder variables, and $0 < S_{\tau_+} \le \xi_{\tau_+}^+$, the result for the nonnegative case implies

$$\begin{aligned}\sup_{b \ge 0} E_\pi R(b) &\le \frac{E_{\pi_+} S_{\tau_+}^2}{E_{\pi_+} S_{\tau_+}} \le \frac{E_{\pi_+} (\xi_{\tau_+}^+)^2}{E_{\pi_+} S_{\tau_+}} \\ &\le \frac{E_{\pi_+} [(\xi_1^+)^2 + \cdots + (\xi_{\tau_+}^+)^2]}{E_{\pi_+} [\xi_1 + \cdots + \xi_{\tau_+}]} = \frac{E_\pi (\xi_1^+)^2}{E_\pi \xi_1}\end{aligned}$$

by Wald's equation for Markov random walks.

When the initial distribution of $X_0$ is $\nu$. Under the assumption of $\sup_x E_x |\xi_1|^2 < \infty$, and $\mu > 0$. It is known that as $b \to \infty$, $E_x(S_{\tau(b)} - b) = E_{\pi_+} S_{\pi_+}^2 / 2 E_{\pi_+} S_{\tau_+} + o(1)$ uniformly in $x \in \mathcal{X}$ [cf. (3.20) of Fuh and Lai (2001)]. Therefore, the difference between $\sup_{b \ge 0} E_\nu R(b)$ and $\sup_{b \ge 0} E_\pi R(b)$ is a constant $K > 0$, and the proof is complete. □

For $z \in \mathbf{C}$ and $\alpha \in \Gamma$, define the operators $\mathbf{P}_z^\alpha, \mathbf{P}, \nu_*^\alpha$ and $\mathbf{Q}$ on $\mathcal{N}$ as (2.11). By (A.1)–(A.3) and Proposition 1 in the Appendix, we can define eigenvalue $\lambda^\alpha(z)$ of the operator $\mathbf{P}_z^\alpha$.

Let "$\mathcal{R}$" and "$\mathcal{I}$" denote "real part of" and "imaginary part of," respectively. Denote $B = [s, s+h)$ and let

$$U^{(\alpha, A)}(B) := \sum_{n=0}^{\infty} P_\nu^\alpha \{ s \le S_n^\alpha \le s + h, X_n^\alpha \in A \} \tag{4.7}$$

be the renewal measure for each $\alpha \in \Gamma$.

LEMMA 3. *Assume the conditions of Theorem 1 hold. For $k \geq 1$, let $\mu_k^\alpha := E_\pi(\xi_1^\alpha)^k > 0$ and $\eta_k^\alpha := \sup_x E_x(\xi_1^\alpha)^k > 0$. Then, for each positive integer $k$, there exist $\mu_{k*}$, $\eta_{k*} > 0$ and $\mu_k^*$, $\eta_k^* < \infty$ such that*

$$\mu_{k*} \leq \inf_{\alpha\in\Gamma} \mu_k^\alpha \leq \sup_{\alpha\in\Gamma} \mu_k^\alpha \leq \mu_k^* \quad \text{and} \quad \eta_{k*} \leq \inf_{\alpha\in\Gamma} \eta_k^\alpha \leq \sup_{\alpha\in\Gamma} \eta_k^\alpha \leq \eta_k^*. \tag{4.8}$$

*Also, there exist $r_1 > 0$ and $\delta > 0$ such that, for all $z$ with $\mathcal{R}(z) \leq r_1$ and $|\mathcal{I}(z)| < \delta$, then for each positive integer $k$, there exists $v_k^* < \infty$ such that, for all $\alpha \in \Gamma$,*

$$\left|\frac{d^k}{dz^k}\lambda^\alpha(z)\right| \leq v_k^*.$$

*Finally, there exists $C$ such that, for all $\alpha \in \Gamma$, and $s \geq 0$ and $h \leq 2$ in (4.7),*

$$U^{(\alpha,A)}(B) \leq C.$$

PROOF. To prove (4.8), we only consider the first part, since the second part can be proved in a similar way. By the assumption of uniformly strong nonlattice in the form (2.3), we have for all $\theta > 0$

$$\widetilde{g}(\theta) := \inf_{\alpha\in\Gamma} |E_\pi e^{i\theta\xi_1^\alpha} - 1| \geq g(\theta) > 0,$$

so that using the fact that $|e^{it} - 1| \leq |t|$ for all real $t$, we obtain for all $\alpha \in \Gamma$ and $\theta > 0$ that

$$\begin{aligned} 0 < \widetilde{g}(\theta) &\leq \int_{x,y\in\mathcal{X}} \int_{[0,\infty)} |e^{i\theta s} - 1| P^\alpha(x, dy \times ds)\pi^\alpha(dx) \\ &\leq \int_{x,y\in\mathcal{X}} \int_{[0,\infty)} \theta s P^\alpha(x, dy \times ds)\pi^\alpha(dx) = \theta\mu_1^\alpha, \end{aligned}$$

where $P^\alpha(\cdot, \cdot \times \cdot)$ denotes the transition probability of $\{(X_n^\alpha, S_n^\alpha), n \geq 0\}$, and $\pi^\alpha(\cdot)$ denotes the invariant probability of $\{(X_n^\alpha, S_n^\alpha), n \geq 0\}$. This implies that $\inf_{\alpha\in A} \mu_1^\alpha \geq \sup_{\theta>0} \widetilde{g}(\theta)/\theta > 0$. Hence, we get the existence of $\mu_{1*}$. The existence of $\mu_{k*}$ for positive integers $k$ now follows from Jensen's inequality.

For the upper bound in (4.8), note that since $e^t > t^k/k!$ for all $t > 0$, we have by K6 that

$$\begin{aligned} \mu_k^\alpha &= \int_{x,y\in\mathcal{X}} \int_{[0,\infty)} s^k P^\alpha(x, dy \times ds)\pi^\alpha(dx) \\ &\leq \frac{k!}{r_1^k} \int_{x,y\in\mathcal{X}} \int_{[0,\infty)} e^{r_1 s} P^\alpha(x, dy \times ds)\pi^\alpha(dx) \\ &\leq \frac{k!C}{r_1^k} \end{aligned}$$

for all $\alpha \in \Gamma$.

To prove the second assertion, note that by Proposition 1 and the Cauchy–Schwarz inequality, for $\mathcal{R}(z) \le r_1/2$ and $|\mathcal{I}(z)| < \delta$, there exists $c > 0$ such that, for all $\alpha \in \Gamma$,

$$\left|\frac{d^k}{dz^k}\lambda^\alpha(z)\right| = |E_\pi^\alpha((\xi_1^\alpha)^k e^{z\xi_1^\alpha})| + c \le E_\pi^\alpha((\xi_1^\alpha)^k e^{r_1\xi_1^\alpha/2}) + c$$
$$\le [E_\pi^\alpha((\xi_1^\alpha)^{2k}) E_\pi^\alpha(e^{r_1\xi_1^\alpha})]^{1/2} + c \le (\mu_{2k}^* C)^{1/2} + c := v_k^*.$$

The final assertion can be proved by using Lemma 1(i), Lemma 2 and (4.8). Let $A = \mathcal{X}$ for simplicity; then for all $\alpha \in \Gamma$, $s \ge 0$ and $h \le 2$, there exists $C_1 < \infty$ such that

$$\begin{aligned}
U^{\alpha,\mathcal{X}}(B) &\le \frac{1}{E_\pi^\alpha \xi_1^\alpha}[E_\nu^\alpha S_{\tau(s+h)}^\alpha - E_\nu^\alpha S_{\tau(s)}^\alpha] + C_1 \\
&= \frac{1}{E_\pi^\alpha \xi_1^\alpha}[h + E_\nu^\alpha R(s+h) - E_\nu^\alpha R(s)] + C_1 \\
&\le \frac{1}{\mu_{1*}}\left[2 + \sup_{\alpha\in\Gamma}\sup_{0\le s<\infty} E_\nu^\alpha R(s)\right] + C_1 \\
&\le \frac{1}{\mu_{1*}}\left[2 + \sup_{\alpha\in\Gamma}\frac{E_\pi(\xi_1^\alpha)^2}{E_\pi \xi_1^\alpha} + K\right] + C_1 \\
&\le \frac{1}{\mu_{1*}}\left[2 + \frac{\mu_2^*}{\mu_{1*}} + K\right] + C_1 := C. \qquad \square
\end{aligned}$$

LEMMA 4. *Assume the conditions of Theorem 1 hold. Then there exist* $r_1 > 0$, $\delta > 0$ *and* $\varepsilon > 0$ *such that, for all* $z$ *satisfying* $0 < \mathcal{R}(z) \le r_1$ *and* $|\mathcal{I}(z)| \le \delta$, *and for all* $\alpha \in \Gamma$, $\lambda^\alpha(z) \ne 1$, *and for all* $z$ *satisfying* $\mathcal{R}(z) = r_1$ *and* $|\mathcal{I}(z)| \le \delta$, *and for all* $\alpha \in \Gamma$,

$$|\lambda^\alpha(z) - 1| \ge \varepsilon. \tag{4.9}$$

PROOF. Let $\mu_k^\alpha := E_\pi(\xi_1^\alpha)^k$. By integration by parts and Proposition 1,

$$\lambda^\alpha(z) = 1 + \mu_1^\alpha z + \int_0^z t\lambda^{\alpha(2)}(z-t)\,dt$$

at least for all $z$ such that $\mathcal{R}(z) < r_1$ and $|\mathcal{I}(z)| \le \delta$, where $^{(2)}$ denotes the second derivative. Therefore, for all $z$ in the set $S := \{z : \mathcal{R}(z) \le r_1/2, |\mathcal{I}(z)| \le \delta, |z| \le \mu_{1*}/v_2^*\}$, where $\mu_{1*}, v_2^*$ are defined in Lemma 3, and all $\alpha \in \Gamma$, we have

$$|\lambda^\alpha(z) - 1| \ge |\mu_1^\alpha z| - \left|\int_0^z t\lambda^{\alpha(2)}(z-t)\,dt\right|$$
$$\ge \mu_{1*}|z| - \frac{v_2^*}{2}|z|^2 \ge \frac{\mu_{1*}}{2}|z|.$$

Take $\varepsilon > 0$ such that the square $S_\varepsilon := \{\mathcal{R}(z) \le \varepsilon, |\mathcal{I}(z)| \le \varepsilon\}$ is contained in the set $S$; for example, $\varepsilon = r_1/2 \wedge \delta \wedge \mu_{1*}/(\sqrt{2}v_2^*)$ will do. Then

$$|\lambda^\alpha(z) - 1| \ge \frac{\mu_{1*}}{2}|z| \qquad \text{for all } z \in S_\varepsilon. \tag{4.10}$$

By the assumption that $\{S_n^\alpha : \alpha \in \Gamma\}$ is uniformly strong nonlattice and Proposition 1, we have, for all $|\theta| < \delta$, $|\lambda^\alpha(i\theta)| = |E_\pi^\alpha e^{i\theta\xi_1^\alpha}| + O(|\theta|)$; this implies that $|\lambda^\alpha(i\theta) - 1| \ge g(\varepsilon) > 0$, for all $|\theta| \ge \varepsilon$ and all $\alpha \in \Gamma$. Take $r := g(\varepsilon)(2v_1^*) \wedge \varepsilon > 0$. Then, for all $0 \le u \le r < r_1$, $\delta > |\theta| \ge \varepsilon$ and $\alpha \in \Gamma$,

$$\begin{aligned}
&|\lambda^\alpha(u + i\theta) - 1| \\
&\quad \ge |\lambda^\alpha(i\theta) - 1| - |\lambda^\alpha(u + i\theta) - \lambda^\alpha(i\theta)| \\
&\quad \ge g(\varepsilon) - \left| \int_{i\theta}^{u+i\theta} \lambda^{\alpha(2)}(z)\,dz \right| \ge g(\varepsilon) - uv_1^* \\
&\quad \ge g(\varepsilon) - \frac{g(\varepsilon)}{2v_1^*}v_1^* = \frac{g(\varepsilon)}{2}.
\end{aligned}$$

Furthermore, (4.10) implies that $|\lambda^\alpha(u + i\theta) - 1|$ is positive for all $0 < u \le r$, $|\theta| \le \varepsilon$ and $\alpha \in \Gamma$, and $|\lambda^\alpha(r + i\theta) - 1|$ is at least $\mu_1^* r/2$ for all $|\theta| \le \varepsilon$. Thus, taking $\delta := g(\varepsilon)/2 \wedge \mu_1^* r/2 > 0$, the lemma is proved. □

Since the rate of convergence in the uniform Markov renewal theorem in Theorem 1 can be proved for each $\alpha \in \Gamma$, the uniformity in $\alpha \in \Gamma$ appealing to Lemmas 3 and 4 when necessary, we will present the proofs by omitting $\alpha$ for simplicity.

Let $B = [s, s + h)$, and recall that $U^{(\alpha,A)}(B)$ defined in (4.7) is the renewal measure. For simplicity, we delete $\alpha$ and denote

$$U^{(A)}(B) := \sum_{n=0}^{\infty} P_\nu\{s \le S_n \le s + h, X_n \in A\} \tag{4.11}$$

as the renewal measure of $\{(X_n, S_n), n \ge 0\}$.

To prove Theorem 1, we evaluate the Fourier transform of the renewal measure $U^A$. As in Carlsson (1983) and Carlsson and Wainger (1982), we perform Fourier inversion of the Fourier transform as a generalized function. We refer the reader to Gelfand and Shilov (1964), Schwartz (1966) and Strichartz (1994) for the basic theory; in particular, the following notation and concepts will be used.

A *test function* $\varphi(s)$ is an infinitely differentiable function that vanishes outside a bounded region in $\mathbf{R}$. Let $\mathcal{D}$ denote the linear space of all test functions, and $\mathcal{D}'$ the space of linear functionals on $\mathcal{D}$. A sequence $\varphi_n \in \mathcal{D}$ is said to converge to zero if $\varphi_n$ and all its derivatives converge to 0 uniformly and vanish outside a common bounded subset of $\mathbf{R}$. A *generalized* *function* is a continuous linear functional on $\mathcal{D}$. A function $f$ defined on $\mathbf{R}$ for which $\int f(s)\varphi(s)\,ds$ is absolutely convergent for any $\varphi\in\mathcal{D}$ is called *locally integrable*. A $C^\infty$ function $f$ on $\mathbf{R}$ is of class $\mathcal{T}$ if $f$ and all its partial derivatives are *rapidly decreasing* in the sense that they are of order $O(|s|^{-a})$ as $|s|\to\infty$, for every $a>0$. Linear functionals on $\mathcal{T}$ are called *tempered distributions*, and $\mathcal{T}'$ denotes the set of all tempered distributions. The Fourier transform $\widehat{\varphi}$ of a function $\varphi\in\mathcal{D}$ is defined by $\widehat{\varphi}(\theta)=\int\varphi(s)\exp(i\theta s)\,ds$ for $\theta\in\mathbf{R}$. The Fourier transform of a generalized function $f$ is the linear functional $\widehat{f}$ defined on the space $\{\psi:\psi$ is the Fourier transform of some $\varphi\in\mathcal{D}\}$ by $(2\pi)(f,\varphi)=(\widehat{f},\widehat{\varphi})$ for all $\varphi\in\mathcal{D}$.

As in renewal theory for simple random walks, the proof of Theorem 1 requires detailed analysis of the characteristic function for the additive component $S_n$. The analysis can be decomposed in two parts: for $|\theta|$ near zero and for $|\theta|$ away from zero. The rate of convergence for the renewal measure to the Lebesgue measure scaled by the mean is given by the analysis of $|\theta|$ near zero. The contribution of $|\theta|$ away from zero is negligible via the property of local integrability. That is, we need to show that for $\theta\in\mathbf{R}$, there exists a $\delta>0$ with $|\theta|>\delta$, $\sum_{n=0}^{\infty}E_\pi(e^{i\theta S_n})$ and its $k$th derivatives with respect to $\theta$ are locally integrable for $k=1,2,\ldots.$ By using $k$th integration by parts, we thus need to show the following lemma.

LEMMA 5. *Suppose* $\{(X_n,S_n),n\ge 0\}$ *is a strongly nonlattice Markov random walk satisfying* (2.4), C1–C4 *and* C6. *Then for every* $c>0$, *for any* $r\ge 2$ *and* $k=0,1,\ldots,r$,

$$\sup_{|\theta|>c}\sum_{n=0}^{\infty}|E_\pi(e^{i\theta S_n}S_n^k)|<\infty. \tag{4.12}$$

PROOF. For $k=0,\ldots,r$ and for any $x\in\mathcal{X}$, we have

$$|E_x(S_n^k e^{i\theta S_n})|\le\sum|E_x(\xi_{j_1}\cdots\xi_{j_k}e^{i\theta S_n})|,$$

where the summation extends over $j_1,\ldots,j_k\in\{1,\ldots,n\}$. There are $n^k$ terms. We shall give upper bounds for each term. Fix $j_1^0,\ldots,j_k^0\in\{1,\ldots,n\}$, and a natural number $m$ to be determined later. Let

$$J=\{j\in\{1,\ldots,n\}:|j-j_p^0|\ge 3m,p=1,\ldots,k\}.$$

Divide $J$ into blocks $A_1,B_1,\ldots,A_l,B_l$ as follows: define $j_1,\ldots,j_l$ by

$$j_1=\inf J\quad\text{and}\quad j_{p+1}=\inf\{j\ge j_p+7m:j\in J\}$$

and let $l$ be the smallest integer for which the inf is undefined. Write

$$A_p=\prod\{e^{i\theta n^{-1/2}\xi_j}:|j-j_p|\le m\},\qquad p=1,\ldots,l,$$

$$B_p = \prod\{e^{i\theta n^{-1/2}\xi_j} : j_p + m + 1 \le j \le j_{p+1} - m - 1\}, \qquad p = 1, \ldots, l-1,$$

$$B_l = \prod\{e^{i\theta n^{-1/2}\xi_j} : j > j_l + m + 1\},$$

$$R = (\xi_{j_1^0} \cdots \xi_{j_k^0}) \prod\{e^{i\theta\xi_j} : j \notin J\}.$$

Then

$$\xi_{j_1^0} \cdots \xi_{j_k^0} e^{i\theta S_n} = \prod_1^l A_p B_p R.$$

We have

$$\begin{aligned}
&\left| E_x R \prod_1^l A_p B_p - E_x R \prod_1^l B_p E(A_p | \xi_j : j \ne j_p) \right| \\
&\quad \le \sum_{q=1}^l \left| E_x R \prod_1^{q-1} A_p B_p (A_q - E(A_p | \xi_j : j \ne j_q)) \right. \\
&\qquad\qquad \left. \times \prod_{q+1}^l B_p E(A_p | \xi_j : j \ne j_p) \right| \qquad (4.13) \\
&\quad \le \sum_{q=1}^l \left| E_x R \prod_1^{q-1} A_p B_p (A_q - E(A_p | \xi_j : j \ne j_q)) \right. \\
&\qquad\qquad \left. \times \prod_{q+1}^l B_p E(A_p | \xi_j : 0 < |j - j_p| \le 3m) \right|.
\end{aligned}$$

The first summation term in (4.13) vanishes since

$$R \prod_1^{q-1} A_p B_p \quad \text{and} \quad \prod_{q+1}^l B_p E(A_p | \xi_j : 0 < |j - j_p| \le 3m)$$

are both measurable with respect to the $\sigma$-field generated by $\xi_j : j \ne j_q$.

Recall that the functions

$$E(A_p | \xi_j : 0 < |j - j_p| \le 3m) \qquad \text{for } p = 1, \ldots, l$$

are weakly dependent since $j_{p+1} - j_p \ge 7m, p = 1, \ldots, l-1$. Using condition C1 we obtain

$$\begin{aligned}
&\left| E_x R \prod_1^l B_p E(A_p | \xi_j : 0 < |j - j_p| \le 3m) \right| \\
&\quad \le (2n^\beta)^r E_x \left| \prod_1^l E(A_p | \xi_j : 0 < |j - j_p| \le 3m) \right|
\end{aligned}$$

$$\leq (2n^{\beta})^r \prod_1^l E_x|E(A_p|\xi_j : 0 < |j - j_p| \leq 3m)| + (2n^{\beta})^r l \cdot 4d^{-1}e^{-dm}.$$

With the strong nonlattice condition (2.3), the conditional strong nonlattice condition (2.4) and Lemma 2 in Statulevicius (1969), we find an upper bound for

$$E_x|E(A_p|\xi_j : 0 < |j - j_p| \leq 3m)|.$$

We have for $|\theta| \geq \delta$, the relation $E_x|E(A_p|\xi_j : j \neq j_q)| \leq e^{-\delta}$, and hence by (2.4), for all $\theta \in \mathbf{R}$, $|\theta| \leq \delta$,

$$E_x|E(A_p|\xi_j : j \neq j_q)| \leq \exp(-\delta|\theta|^2/n).$$

Therefore, for all $\theta \in \mathbf{R}$,

$$\begin{aligned} E_x|E(A_p|\xi_j : 0 < |j - j_p| \leq 3m)| &\leq E_x|E(A_p|\xi_j : j \neq j_q)| \\ &\leq \max(\exp(-\delta|\theta|^2/n), e^{-\delta}). \end{aligned}$$

If we choose $K$ appropriately and let $m$ be the integral part of $K \log n$, then the assertion of the lemma follows from

$$\begin{aligned} \exp(-\delta|\theta|^2/n)^{n/m} &\leq \exp(-\delta|\theta|^2/(K \log n)) \\ &\leq \exp(-\delta' n^{\varepsilon/2}) \end{aligned}$$

for $|\theta| \geq cn^{\varepsilon}$ and some $\delta' > 0$. □

Proof of Theorem 1. By using the same argument as that in Theorem 2.4 of Fuh and Lai (2001), we have (2.5) and (2.6). The details are omitted.

To prove (2.7), let $g(s) \in \mathbf{C}^{\infty}$ and have support in $\Omega = \{s : |s| \leq 1\}$. Let $g_{\varepsilon}(s) = \varepsilon g(s/\varepsilon)$, let $I_A$ be the indicator function of the set $A$ and let $\Omega_c = \{s : |s| \leq c\}$ for $c > 0$. Let $L$ be the measure with density $ds/\mu$. As in Carlsson and Wainger (1982), we have that

$$\begin{aligned} g_{\varepsilon} * I_{\Omega_{1-\varepsilon}} * (U^{(A)} - L)(s) - K\varepsilon &\leq (U^{(A)} - L)(\Omega + s) \\ &\leq g_{\varepsilon} * I_{\Omega_{1+\varepsilon}} * (U^{(A)} - L)(s) + K\varepsilon, \end{aligned} \tag{4.14}$$

where $K$ can be chosen uniformly for $c$ bounded. Letting $\varepsilon = s^{-k}$ with $k$ large enough, we have (2.7) if we can show that there exists $r > 0$ such that

$$|g_{\varepsilon} * I_{\Omega_c} * (U^{(A)} - L)(s)| \leq Ke^{-rs}. \tag{4.15}$$

To prove (4.15), we consider the Fourier transform

$$\psi(\theta) = \sum_{n=0}^{\infty} E_\nu\{e^{i\theta S_n} h_1(X_n)\} \tag{4.16}$$
$$= \sum_{n=0}^{\infty} \lambda^n(\theta)\nu_* \mathbf{Q}_\theta h_1 + \sum_{n=0}^{\infty} \nu_* \mathbf{P}_\theta^n (I - \mathbf{Q}_\theta) h_1,$$

where $\mathbf{P}_\theta, \mathbf{Q}_\theta$ are defined in (2.11) with $z = i\theta$, and $h_1 := I_A$. Note that the second equation in (4.16) follows from Proposition 1(i).

As in Carlsson and Wainger (1982), there exists a $\delta > 0$ such that, for $|\theta| < \delta$,

$$\psi(\theta) = \left(\frac{1}{1-\lambda(\theta)} + \frac{\pi}{\mu}\delta(\theta)\right)\nu_* \mathbf{Q}_\theta h_1 + \sum_{n=0}^{\infty} \nu_* \mathbf{P}_\theta^n (I - \mathbf{Q}_\theta) h_1, \tag{4.17}$$

where $\delta(\cdot)$ denotes the Dirac delta function. By using Fourier inversion of the generalized function, we have

$$g_\varepsilon * I_{\Omega_c} * U^{(A)}(s)$$
$$= \frac{1}{2\pi}\int e^{-i\theta s} \hat{g}_\varepsilon(\theta) \hat{I}_{\Omega_c}(\theta) \psi(\theta)\, d\theta$$
$$= \frac{\pi}{\mu} g(0) + \frac{1}{2\pi}\int_{0<|\theta|<\delta} e^{-i\theta s} \frac{\hat{g}_\varepsilon(\theta)\hat{I}_{\Omega_c}(\theta)}{1-\lambda(\theta)} \nu_* \mathbf{Q}_\theta h_1 \, d\theta \tag{4.18}$$
$$+ \frac{1}{2\pi}\int_{0<|\theta|<\delta} e^{-i\theta s} \hat{g}_\varepsilon(\theta)\hat{I}_{\Omega_c}(\theta) \sum_{n=0}^{\infty} \nu_* \mathbf{P}_\theta^n (I - \mathbf{Q}_\theta) h_1 \, d\theta \tag{4.19}$$
$$+ \frac{1}{2\pi}\int_{|\theta|>\delta} e^{-i\theta s} \hat{g}_\varepsilon(\theta)\hat{I}_{\Omega_c}(\theta) \sum_{n=0}^{\infty} E_\nu\{e^{i\theta S_n} h_1(X_1)\}\, d\theta, \tag{4.20}$$

where $\hat{g}_\varepsilon$ denotes the Fourier transform of $g_\varepsilon$ and $\hat{I}_{\Omega_c}$ denotes the Fourier transform of $I_{\Omega_c}$.

Equation (4.18) can be analyzed by making use of the Taylor expansion of $\lambda(\theta)$ as that in Proposition 1, which says that for any $r \geq 2$ as $|\theta| \to 0$,

$$\left(\frac{1}{|\theta|}(\lambda(\theta) - 1)\right)^{(k)} = O(1) \qquad \text{for } k \leq r-1,$$
$$\left(\frac{1}{|\theta|}(\lambda(\theta) - 1)\right)^{(r)} = o\left(\frac{1}{|\theta|}\right),$$
$$\left(\frac{1}{|\theta|^2}(\lambda(\theta) - 1 + i\theta\mu)\right)^{(k)} = O(1) \qquad \text{for } k \leq r-2,$$
$$\left(\frac{1}{|\theta|^2}(\lambda(\theta) - 1 + i\theta\mu)\right)^{(k)} = o(|\theta|^{r-k-2}) \qquad \text{for } k = r-1, r,$$

where ${}^{(k)}$ denotes the $k$th derivative. Also $\nu_* \mathbf{Q}_\theta h_1$ and its $k$th derivatives converge to 0 as $|\theta| \to 0$.

Next, we want to verify that the rate of convergence in (4.18) is $O(e^{-rs})$ for some $r > 0$ as $s \to \infty$. Note that

$$\frac{1}{2\pi} \int_{0<|\theta|<\delta} \mathcal{R}\{e^{-i\theta s} \hat{g}_\varepsilon(\theta) \hat{I}_{\Omega_c}(\theta) \nu_* \mathbf{Q}_\theta h_1 (i/\mu\theta)\} \, d\theta = \frac{g(0)}{2\mu} + O(e^{-rs}), \tag{4.21}$$

as $s \to \infty$. Consider

$$\begin{aligned} &\frac{1}{2\pi} \int_{0<|\theta|<\delta} e^{-i\theta s} \hat{g}_\varepsilon(\theta) \hat{I}_{\Omega_c}(\theta) \nu_* \mathbf{Q}_\theta h_1 \left( \frac{1}{1-\lambda(i\theta)} - \frac{i}{\mu\theta} \right) d\theta \\ &\quad = \frac{1}{2\pi} \lim_{\varepsilon \to 0} \biggl( \int_{-\delta<\theta<\varepsilon} e^{-i\theta s} \hat{g}_\varepsilon(\theta) \hat{I}_{\Omega_c}(\theta) \nu_* \mathbf{Q}_\theta h_1 \\ &\qquad \times \left( \frac{1}{1-\lambda(i\theta)} - \frac{i}{\mu\theta} \right) d\theta \\ &\qquad + \int_{\varepsilon<\theta<\delta} e^{-i\theta s} \hat{g}_\varepsilon(\theta) \hat{I}_{\Omega_c}(\theta) \nu_* \mathbf{Q}_\theta h_1 \\ &\qquad \times \left( \frac{1}{1-\lambda(i\theta)} - \frac{i}{\mu\theta} \right) d\theta \biggr). \end{aligned} \tag{4.22}$$

Let $0 < u_1 \in \Theta$, where $\Theta$ is defined in K6. For any $u \in (0, u_1)$, and $z = u + i\theta \in \mathbf{C}$, consider four lines $L_1(\varepsilon) = \{z : \mathcal{R}(z) = 0, \varepsilon \le |\mathcal{I}(z)| \le \delta\}$, $L_2 = \{z : \mathcal{R}(z) \in [0, u], \mathcal{I}(z) = \delta\}$, $L_3 = \{z : \mathcal{R}(z) \in [0, u], \mathcal{I}(z) = -\delta\}$, $L_4 = \{z : \mathcal{R}(z) = u, |\mathcal{I}(z)| \le \delta\}$, and one semicircle, $L_5(\varepsilon)$, from $\varepsilon$ to $-\varepsilon$, oriented clockwise. Define

$$h(z) = e^{-zs} \hat{g}_\varepsilon(z) \hat{I}_{\Omega_c}(z) \nu_* \mathbf{Q}_z h_1 \left( \frac{1}{1-\lambda(z)} - \frac{i}{\mu z} \right). \tag{4.23}$$

Since $h$ is analytic in the regions from $L_1(\varepsilon)$ to $L_5(\varepsilon)$, by Cauchy's theorem for contour integral, we have

$$\begin{aligned} &\int_{L_1(\varepsilon)} h(z) \, dz + \int_{L_2} h(z) \, dz + \int_{L_3} h(z) \, dz \\ &\qquad + \int_{L_4} h(z) \, dz + \int_{L_5(\varepsilon)} h(z) \, dz = 0. \end{aligned} \tag{4.24}$$

The continuity $h$ yields that the residual of $h(z)$ at 0 is 0, whence

$$\lim_{\varepsilon \to 0} \int_{L_5(\varepsilon)} h(z) \, dz = 0.$$

Combining with the Riemann–Lebesgue lemma, we have

$$\lim_{\varepsilon \to 0} \int_{L_1(\varepsilon)} h(z) \, dz$$

$$
\begin{aligned}
&= \int_{-\delta}^{\delta} e^{-(u+i\theta)s} \hat{g}_\varepsilon(u+i\theta) \hat{I}_{\Omega_c}(u+i\theta) \nu_* \mathbf{Q}_{u+i\theta} h_1 \\
&\qquad \times \left( \frac{1}{1-\lambda(u+i\theta)} - \frac{i}{\mu(u+i\theta)} \right) d\theta \qquad (4.25)\\
&= O(e^{-rs}) \qquad \text{for some } \ r>0 \text{ as } s\to\infty.
\end{aligned}
$$

To analyze (4.20), we make use of Lemma 5 which implies that $\sum_{n=0}^{\infty} E_\nu\{e^{i\theta S_n}\}$ and its partial derivatives up to order $r$ are bounded for any $r \geq 2$, and a fortiori locally integrable, in the region $\{\theta : |\theta| \geq \delta\}$. Moreover, $\hat{I}_{\Omega_c}$ and its derivatives are bounded by a constant times $\prod_{i=1}^{d} |\theta|^{-1}$ as $|\theta| \to \infty$. Therefore $N$ integrations by parts as in page 359 of Carlsson and Wainger (1982) can be used to show that $(4.20) = O(|\log \varepsilon|^d / s_1^N)$ for any $N$. Since $\eta(\theta)$ and its partial derivatives of order $r$ are bounded for $0 < |\theta| \leq \delta^*$ for any $r \geq 2$, $N$ integrations by parts can also be used to show that $(4.19) = O(s_1^{-N})$ for any $N$. Therefore, by making use (4.18)–(4.20), (4.23) and (4.25), we have (4.15) and hence get the proof. $\square$

**5. Proof of Theorem 2.** In this section, we assume the Markov random walk $\{(X_n, S_n), n \geq 0\}$ defined as (1.1) is uniformly ergodic with respect to a given norm, and the stationary mean $\mu = 0$. Under the minorization condition C1, making use of the results in Hipp (1985), Malinovskii (1987) and Jensen (1989), we have the following asymptotic expansions of the density for the distribution in Markov random walks.

LEMMA 6. *Assume* C1–C5 *with* $r = 3$ *in* C3. *We assume, without loss of generality, the asymptotic variance* $\sigma^2 = 1$. *Then*

$$
P_\nu\{S_n \leq s\sqrt{n}\} = \Phi(s) + \phi(s)\frac{Q_1(s)}{\sqrt{n}} + (1+|s|^3)^{-1} o\left(\frac{1}{\sqrt{n}}\right), \tag{5.1}
$$

*where* $o(\cdot)$ *is uniform in* $s$. *Here* $\Phi(\cdot)$ *denotes the standard normal distribution,* $\phi(\cdot)$ *denotes the standard normal density, and* $Q_1(s) = \kappa/6(1-s^2) + \kappa_\nu$, *where* $\kappa = E_\pi \xi_1^3 + 3\sum_{t=1}^{\infty} E_\pi \xi_1^2 \xi_{t+1} + 3\sum_{t=1}^{\infty} E_\pi \xi_1 \xi_{t+1}^2 + 6\sum_{t_1,t_2=1}^{\infty} E_\pi \xi_1 \xi_{t_1+1} \xi_{t_1+t_2+1}$ *and* $\kappa_\nu = \sum_{t=1}^{\infty} E_\nu \xi_t$. *Note that* $k_\nu = 0$ *if* $\nu = \pi$.

*Furthermore, if* $P_\pi\{S_n \leq s\sqrt{n}\}$ *has a density* $p_{\pi,n}(s\sqrt{n})$, *then*

$$
p_{\pi,n}(s\sqrt{n})\sqrt{n} = \phi(s)\left(1 + \frac{k}{6\sqrt{n}}(s^3 - 3s)\right) + (1+|s|^3)^{-1} o\left(\frac{1}{\sqrt{n}}\right), \tag{5.2}
$$

*where* $o(\cdot)$ *is uniform in* $s$.

In the following, we shall assume C1–C5 and C7 hold. Lemma 7 is taken from Theorem 5 of Fuh and Lai (1998); we include it here for completeness.

LEMMA 7. *Let* $r \geq 1$. *Assume* $\sup_x E_x(\xi_1^+)^{r+1} < \infty$, *where* $\xi_1^+$ *denotes the positive part of* $\xi_1$. *Furthermore, assume there exists* $\varepsilon > 0$ *such that* $\inf_x P_\pi\{\xi_1 \leq -\varepsilon | X_1 = x\} > 0$. *Then,* $E_\pi S_{\tau_+}^r < \infty$ *and* $\{S_{\tau(b)} - b, b > 0\}$ *is uniformly integrable under the probability* $P_\pi$.

A Markov random walk is called *lattice* with span $d > 0$ if $d$ is the maximal number for which there exists a measurable function $\gamma : \mathcal{X} \to [0, \infty)$, called the shift function, such that $P\{\xi_1 - \gamma(x) + \gamma(y) \in \{\dots, -2d, -d, 0, d, 2d, \dots\} | X_0 = x, X_1 = y\} = 1$ for almost all $x, y \in \mathcal{X}$. If no such $d$ exists, it is called *nonlattice*.

Let $W(t)$, $0 \leq t \leq \infty$, denote Brownian motion with drift $\mu$ and put $\tau_W = \tau_W(b) = \inf\{t : W(t) \geq b\}$. Define the inverse Gaussian distribution $G(t; \mu, b) = P^{(\mu)}\{\tau_W(b) \leq t\}$ and $H_+(s) = (E_{\pi_+} S_{\tau_+})^{-1} \int_0^s P_{\pi_+}\{S_{\tau_+} \geq t\}\, dt$.

LEMMA 8. *Assume* $P$ *is nonlattice. Suppose* $b \to \infty$ *and* $m \to \infty$ *so that, for some fixed* $0 < \zeta < \infty$, $b = \zeta m^{1/2}$. *Then for all* $0 \leq t, s \leq \infty$,

$$P_\pi\{\tau(b) \leq mt,\ S_{\tau(b)} - b \leq s\} \longrightarrow G(t; 0, \zeta) H_+(s).$$

PROOF. Note that

$$\begin{aligned} &P_\pi\{\tau(b) > mt,\ S_{\tau(b)} - b \leq s\} \\ &\quad = E_\pi(P_\pi\{S_{\tau(b)} - b \leq s | \tau(b) > mt, S_m\}; \tau(b) > mt). \end{aligned}$$

Let $\tilde{b} = b - b^{1/4}$. Then by the central limit theorem for Markov random walks [cf. Theorem 17.2.2 of Meyn and Tweedie (1993)], we have

$$\begin{aligned} &E_\pi(P_\pi\{S_{\tau(b)} - b \leq s | \tau(b) > mt, S_m\}; \tau(b) > mt, \tilde{b} < S_m < b) \\ &\quad \leq P_\pi\{\tilde{b} < S_m < b\} \to 0. \end{aligned}$$

Moreover, it is also known that $P_{\pi_+}\{\tau_+ < \infty\} = 1$ and by Lemma 7, $E_{\pi_+} S_{\tau_+} < \infty$, and

$$\begin{aligned} &P_\pi\{S_{\tau(b)} - b \leq s | \tau(b) > mt, S_m \leq \tilde{b}\} \\ &\quad = \int_0^{\tilde{b}} \int_{\mathcal{X}} P_\pi\{S_{\tau(b)} - b \leq s | X_m = x, \tau(b) > mt, S_m \in dv\} \\ &\qquad\quad \times P_\pi\{X_m \in dx | \tau(b) > mt, S_m \in dv\} \\ &\quad = \int_0^{\tilde{b}} \int_{\mathcal{X}} P_x\{S_{\tau(b-v)} - (b - v) \leq s\} P_\pi\{X_m \in dx | \tau(b) > mt, S_m \in dv\}. \end{aligned}$$

Since $P_\pi\{\tau(b) > mt, S_m \leq \tilde{b}\} > 0$ as $b \to \infty$, therefore for $v$ uniformly in $\{\tau(b) > mt, S_m \leq \tilde{b}\}$, $P_x\{S_{\tau(b-v)} - (b - v) \leq s\} \to H_+(s)$ as $b \to \infty$ by the Markov renewal theorem in Theorem 1. Hence, uniformly in $\{\tau(b) > mt, S_m \le \tilde{b}\}$, as $b \to \infty$,

$$P_\pi\{S_{\tau(b)} - b \le s | \tau(b) > mt, S_m\} \to H_+(s).$$

Under irreducible assumption and the minorization condition C1, for $r = 2$, assumption C3 ensures that the Poisson equation (4.1) has a solution $\Delta$ which satisfies $E_\pi(\Delta(X_1))^2 < \infty$. Therefore, by the invariance principle for Harris recurrent Markov chain [cf. Theorem 17.4.4 of Meyn and Tweedie (1993)], the limiting marginal distribution of $\tau(b)/m$ is $G(t; 0, \zeta)$. Hence,

$$\begin{aligned} &P_\pi\{\tau(b) > mt,\ S_{\tau(b)} - b \le s\} \\ &\quad = E_\pi(P_\pi\{S_{\tau(b)} - b \le s | \tau(b) > mt, S_m\}; \tau(b) > mt, S_m \le \tilde{b}) + o(1) \\ &\quad = H_+(s) P_\pi\{\tau(b) > mt, S_m \le \tilde{b}\} + o(1) \to H_+(s)(1 - G(t; 0, \zeta)). \quad \square \end{aligned}$$

Lemma 9. *Let $0 < \zeta < \infty$, and $R_m = S_\tau(b) - \zeta m^{1/2}$. Then for any $\varepsilon > 0$, $P_\pi\{R_m > \varepsilon m^{1/2}\} = o(m^{-1/2})$.*

Proof. By Markov's inequality we have

$$P_\pi\{R_m > \varepsilon m^{1/2}\} = \varepsilon^{-1} m^{-1/2} \int_{\{R_m > \varepsilon m^{1/2}\}} R_m \, dP_\pi$$

which is $o(m^{-1/2})$ by Lemma 7. $\square$

Lemma 10. *Let $0 < \zeta < \infty$, $0 \le s \le \infty$, and $m_1 = m(1 - (\log m)^{-2})$. Then as $m \to \infty$,*

$$P_\pi\{m_1 < \tau < m, S_m < (\zeta - s) m^{1/2}\} = o(m^{-1/2}), \tag{5.3}$$

$$P_\pi\{m_1 < \tau \le m, S_m \ge (\zeta + s) m^{1/2}\} = o(m^{-1/2}). \tag{5.4}$$

Proof. By Lemma 9,

$$\begin{aligned} &P_\pi\{m_1 < \tau \le m, S_m \ge (\zeta + s) m^{1/2}\} \\ &\quad = P_\pi\{m_1 < \tau \le m, R_m < \tfrac{1}{2} s m^{1/2}, S_m \ge (\zeta + s) m^{1/2}\} + o(m^{-1/2}) \qquad (5.5) \\ &\quad \le \sup_{m_1 < n \le m} P_\pi\{S_{m-n} \ge \tfrac{1}{2} s m^{1/2}\} + o(m^{-1/2}). \end{aligned}$$

It is easy to see (5.5) is $o(m^{-1/2})$. By using a similar argument, we have that (5.3) is $o(m^{-1/2})$. $\square$

Proof of Theorem 2. Since we will consider time delay in the proof, we denote $S_0 = s_0$ for convenience. Let $P_\pi^{(m,s_0,s)}(A) = P_\pi\{A | S_0 = s_0, S_m =$ $s\}$. Set $s_0 = m^{1/2}\lambda_0$ and $s = m^{1/2}\zeta_0$ for some $\lambda_0, \zeta_0 < \zeta$. Let $s' = 2b - s = m^{1/2}(2\zeta - \zeta_0)$ denote $s$ reflected about $b$. From (5.2) in Lemma 6 and the Markov renewal theorem in Theorem 1, in which $\Gamma$ has only one element, it follows as in Lemmas 7–10 that, for $m_1 = m(1 - (\log m)^{-2})$ and some $\varepsilon_m \to 0$,

$$P_\pi^{(m,s_0,s)}\{\tau < m\} - P_\pi^{(m,s_0,s)}\{\tau < m_1, S_\tau - b < m^{1/2}\varepsilon_m\} = o(m^{-1/2}) \tag{5.6}$$

and

$$P_\pi^{(m,s_0,s')}\{\tau < m\} - P_\pi^{(m,s_0,s')}\{\tau < m_1, S_\tau - b < m^{1/2}\varepsilon_m\} = o(m^{-1/2}). \tag{5.7}$$

Set $A_m = \{\tau < m_1, S_\tau - b < m^{1/2}\varepsilon_m\}$, and let $L^{(m)}(n, S_n)$ denote the likelihood ratio of $\xi_1, \ldots, \xi_n$ under $P_\pi^{(m,s_0,s)}$ relative to $P_\pi^{(m,s_0,s')}$. For all $n \le m - n_0$,

$$L^{(m)}(n, S_n) = \frac{p_{\pi,m-n}(s - S_n)p_{\pi,m}(s' - s_0)}{p_{\pi,m-n}(s' - S_n)p_{\pi,m}(s - s_0)}, \tag{5.8}$$

where $p_{\pi,m}(s)$ is defined in (5.2). By (5.6) and Wald's likelihood ratio identity for Markov chains

$$P_\pi^{(m,s_0,s)}\{\tau < m\} = E_\pi^{(m,s_0,s')}\{L^{(m)}(\tau, S_\tau); A_m\} + o(m^{-1/2}). \tag{5.9}$$

Substitution of (5.8) into (5.9) and expansion with the aid of (5.2) gives the first-order result,

$$P_\pi^{(m,s_0,s)}\{\tau < m\} \to \exp[-2(\zeta - \lambda_0)(\zeta - \zeta_0)].$$

This motivates the following reformulation of (5.9), which is justified by (5.7) and the fact that $P_\pi^{(m,s_0,s')}\{\tau = m\} = o(m^{-1/2})$:

$$\begin{aligned} &P_\pi^{(m,s_0,s)}\{\tau < m\} - \exp[-2(\zeta - \lambda_0)(\zeta - \zeta_0)] + o(m^{-1/2}) \\ &\quad = E_\pi^{(m,s_0,s')}\{L^{(m)}(\tau, S_\tau) - \exp[-2(\zeta - \lambda_0)(\zeta - \zeta_0)]; A_m\}. \end{aligned} \tag{5.10}$$

The likelihood ratio of $\xi_1, \ldots, \xi_n$ under $P_\pi^{(m,s_0,s')}$ relative to $P_\pi\{\cdot|S_0 = s_0\}$ is $p_{\pi,m-n}(s' - S_n)/p_{\pi,m}(s' - s_0)$. Hence by (5.8) and Wald's likelihood ratio identity once again, the right-hand side of (5.10) becomes

$$\begin{aligned} E_\pi\biggl\{&\frac{p_{\pi,m-\tau}(s - S_\tau)}{p_{\pi,m}(s - s_0)} \\ &- \exp[-2(\zeta - \lambda_0)(\zeta - \zeta_0)]\frac{p_{\pi,m-\tau}(s' - S_\tau)}{p_{\pi,m}(s' - s_0)}; A_m \bigg| S_0 = s_0\biggr\}. \end{aligned} \tag{5.11}$$

The rest of the proof of (2.9) involves use of Lemma 6 to expand the integrand in (5.11) and application of Lemma 8 to evaluate the resulting expectation. Let $R_m = S_\tau - m^{1/2}\zeta$. Some tedious algebra gives that the integrand in (5.11) equals

$$[(1-\tau/m)^{1/2}\phi(\zeta_0-\lambda_0)]^{-1} \times \left\{\phi\left(\frac{\zeta-\zeta_0+R_m/m^{1/2}}{(1-\tau/m)^{1/2}}\right)-\phi\left(\frac{\zeta-\zeta_0-R_m/m^{1/2}}{(1-\tau/m)^{1/2}}\right)\right\},$$

which can be expanded to give

$$\begin{aligned} &-2(1-\tau/m)^{-1/2}\exp[\tfrac{1}{2}(\zeta_0-\lambda_0)^2-\tfrac{1}{2}(\zeta-\zeta_0)^2/(1-\tau/m)] \\ &\quad\times[(\zeta-\zeta_0)R_m/m^{1/2}(1-\tau/m)]+o((1+R_m^2)/m) \end{aligned} \tag{5.12}$$

uniformly on $A_m$. According to Lemma 8, $\tau/m$ and $R_m$ are asymptotically independent, converge in law, and by Lemma 7, $R_m$ is uniformly integrable. Also, (5.12) is a bounded, continuous function of $\tau/m$ on $A_m$. Hence, (5.12) can be substituted into (5.11) and Lemma 8 applied to evaluate the result. Putting $s_0=0$ and performing the appropriate integrations yields (2.9).

Formally, (2.10) follows by substituting (2.9) into

$$P_\pi\{\tau<m,S_m<c\}=\int_{(-\infty,c)}P_\pi^{(m,0,s)}\{\tau<m\}P_\pi\{S_m\in ds\}. \tag{5.13}$$

However, some care is required to justify this calculation, especially in the case $c=b$ ($\gamma=\zeta$), when $s$ in (5.13) can be arbitrarily close to $b$. It is easy to see that (2.9) holds uniformly on each compact subinterval of $(-\infty,\zeta)$; but if $\zeta_0\to\zeta$, (5.12) is not necessarily bounded, (5.6) may fail to hold, and indeed the proof of (2.9) disintegrates.

To circumvent this difficulty, we need to apply the duality argument to a time-reversed Markov chain. By condition C7, recall that $P_x(A)=\int_A p(x,y)M(dy)$ for all $A\in\mathcal{A}$, where $p(x,\cdot)=dP_x/dM$. Letting $Q^{x,y}(B)=P(\xi_1\in B|X_0=x,X_1=y)$, we can express the transition probability function (1.1) as

$$P(x,A\times B)=\int_A p(x,y)Q^{x,y}(B)M(dy):=\int_A F(x,y;B)M(dy), \tag{5.14}$$

where $F(x,y;B)=p(x,y)Q^{x,y}(B)$. For ease of notation, we still denote $\pi$ as the density of the invariant probability measure. We shall use $\sim$ to refer to the time-reversed (or dual) chain $\{(\tilde{X}_n,\tilde{S}_n),n\ge 0\}$ with transition kernel:

$$\tilde{F}(y,x;B)=F(x,y;B)\pi(x)/\pi(y). \tag{5.15}$$

Let $\tau^*=\sup\{n:n<m,S_m\ge b\}$, and let $\tilde{\tau}=\inf\{n:\tilde{S}_n\ge b\}$ denote the first passage time of the time-reversed Markov random walk $\tilde{S}_n$ at the linear boundary $b$. Observe that $P_\pi^{(m,0,\zeta)}\{\tau<m\}=P_\pi^{(m,0,\zeta)}\{\tau^*>0\}=P_\pi^{(m,\zeta,0)}\{\tilde{\tau}<$

$m\}$; so to approximate $P_\pi^{(m,0,\zeta)}\{\tau < m\}$ for $\zeta - \varepsilon \le \zeta_0 < \zeta$, it suffices to consider $P_\pi^{(m,\lambda,\zeta)}\{\tilde\tau < m\}$ for $\zeta = 0$ and $\zeta - \varepsilon \le \lambda_0 < \zeta$ (recall that $\lambda_0 = m^{-1/2}s_0$). It is easy to see from (5.10)–(5.12) that uniformly for $\zeta - \varepsilon \le \lambda_0 \le s_0 - m^{-1/2} < \zeta$, $P_\pi^{(m,s_0,0)}\{\tilde\tau < m\} = \exp\{-2\zeta(\zeta - \lambda_0)\} + o(m^{-1/2})$, which suffices to justify formal substitution of (2.9) into (5.13) for $\zeta$ in a neighborhood of $b$ and complete the proof when $\gamma = \zeta$. $\square$

**6. Proof of Theorem 3.** The way to prove Theorem 3 is a suitable application of Theorem 1 via the following lemmas.

LEMMA 11. *Assume the conditions of Theorem* 3 *hold. Then, there exists* $\delta > 0$ *and* $|\alpha| \le \delta$ *such that the induced Markov chain* $\{(X_n^\alpha, S_n^\alpha), n \ge 0\}$ *with transition probability* (2.13) *is aperiodic and irreducible. Moreover, it is uniform ergodic with respect to a given norm and satisfying* K1–K6.

PROOF. For $|\alpha| \le \delta$, it is known [cf. Ney and Nummelin (1987)] that $\{(X_n^\alpha, S_n^\alpha), n \ge 0\}$ is an aperiodic and irreducible Markov chain. Since $\{(X_n, S_n), n \ge 0\}$ satisfies C1, $P^\alpha$ is geometrically $e^{-\Lambda(\alpha)}$-recurrent for $|\alpha| < \delta$, and therefore is $e^{-\Lambda(\alpha)}$-uniformly ergodic, compare Theorem 4.1 of Ney and Nummelin (1987).

Define

$$\Psi^\alpha(dy \times ds) = \frac{\Psi(dy \times ds)e^{-\Lambda(\alpha)+\alpha s}r(y;\alpha)}{(\Psi r)(\alpha)},$$

where $(\Psi r)(\alpha)$ is a normalizing constant, and

$$h^\alpha(x) = (\Psi r)(\alpha)r^{-1}(x;\alpha)h(x),$$

where $\Psi(\cdot)$ and $h(\cdot)$ are defined in C1 with

$$P(x, dy \times ds) \ge h(x)\Psi(dy \times ds).$$

This implies that

$$\begin{aligned} P^\alpha(x, dy \times ds) &\ge \frac{\Psi(dy \times ds)e^{-\Lambda(\alpha)+\alpha s}r(y;\alpha)}{(\Psi r)(\alpha)}(\Psi r)(\alpha)r^{-1}(x;\alpha)h(x) \\ &= h^\alpha(x)\Psi^\alpha(dy \times ds). \end{aligned}$$

Therefore the mixing condition K1 hold.

To prove the moment condition K6, denote $\lambda^\alpha(z)$ as the eigenvalue of $\mathbf{P}_z^\alpha$. By (2.13), we have

$$\begin{aligned} &|E_\pi^\alpha(e^{\theta\xi_1^\alpha}) - E_\pi(e^{\theta\xi_1})| \\ &\quad = \left|\int_{x,y\in\mathcal{X}}\int_{-\infty}^{\infty}\frac{r(y;\alpha)}{r(x;\alpha)}e^{\theta s}\left(\frac{r(x;\alpha)}{r(y;\alpha)} - e^{\alpha s-\Lambda(\alpha)}\right)P^\alpha(x, dy \times ds)\pi^\alpha(dx)\right| \\ &\quad \le \int_{x,y\in\mathcal{X}}\int_{-\infty}^{\infty}\left|\frac{r(y;\alpha)}{r(x;\alpha)}\right||e^{\theta s}|\left|\frac{r(x;\alpha)}{r(y;\alpha)} - e^{\alpha s-\Lambda(\alpha)}\right|P^\alpha(x, dy \times ds)\pi^\alpha(dx) \end{aligned}$$

for $\theta \in \Theta \subset \mathbf{R}$. From $\lim_{\alpha\downarrow 0} |r(y;\alpha)/r(x;\alpha)| = 1$, we get $\lim_{\alpha\downarrow 0} \sup_\theta |E^\alpha_\pi(e^{\theta\xi^\alpha_1}) - E^\alpha_\pi(e^{\theta\xi_1})| = 0$ by dominated convergence. Therefore, we may choose $\alpha^* > 0$ so that $\sup_{\alpha\in[0,\alpha^*]} \sup_\theta |E_\pi(e^{\theta\xi^\alpha_1}) - E_\pi(e^{\theta\xi_1})| \le C$, say. Hence, condition C6 implies K6 holds. By using a similar argument, we also have the moment conditions K2–K4. $\square$

LEMMA 12. *Assume the conditions of Theorem* 3 *hold. Then there exists* $\alpha^* > 0$ *such that the family* $\{(X^\alpha_n, S^\alpha_n), n \ge 0 : 0 \le \alpha \le \alpha^*\}$ *satisfies* (2.3) *and* (2.4).

PROOF. Since the proofs of (2.3) and (2.4) are similar, we only prove (2.3), the uniformly strong nonlattice case. Let $\lambda^\alpha(z)$ denote the eigenvalue of $\mathbf{P}^\alpha_z$. Since $\{(X_n, S_n), n \ge 0\}$ is assumed to be strongly nonlattice, we have that $g(1) := \inf_{|\theta|\ge 1} |1 - E_\pi(e^{i\theta\xi_1})| > 0$. However, by (2.13),

$$\begin{aligned}
&|E_\pi(e^{i\theta\xi_1}) - E^\alpha_\pi(e^{i\theta\xi^\alpha_1})| \\
&\quad = \left| \int_{x,y\in\mathcal{X}} \int_{-\infty}^{\infty} \frac{r(y;\alpha)}{r(x;\alpha)} e^{i\theta s} \left( e^{\alpha s - \Lambda(\alpha)} - \frac{r(x;\alpha)}{r(y;\alpha)} \right) P^\alpha(x, dy \times ds) \pi^\alpha(dx) \right| \\
&\quad \le \int_{x,y\in\mathcal{X}} \int_{-\infty}^{\infty} \left| \frac{r(y;\alpha)}{r(x;\alpha)} \right| \left| e^{\alpha s - \Lambda(\alpha)} - \frac{r(x;\alpha)}{r(y;\alpha)} \right| P^\alpha(x, dy \times ds) \pi^\alpha(dx)
\end{aligned}$$

for all real $\theta$. Since $\lim_{\alpha\downarrow 0} |r(y;\alpha)/r(x;\alpha)| = 1$, so $\lim_{\alpha\downarrow 0} \sup_\theta |E_\pi(e^{i\theta\xi_1}) - E^\alpha_\pi(e^{i\theta\xi^\alpha_1})| = 0$ by dominated convergence theorem. Therefore, we may choose $\alpha^* > 0$ so that $\sup_{\alpha\in[0,\alpha^*]} \sup_\theta |E_\pi(e^{i\theta\xi_1}) - E^\alpha_\pi(e^{i\theta\xi^\alpha_1})| \le g(1)/2$, say. Choosing $\alpha^*$ in this way, by applying the definition of $g(1)$ and the triangle inequality, we obtain

$$\inf_{\alpha\in[0,\alpha^*]} \inf_{|\theta|\ge 1} |1 - E^\alpha_\pi(e^{i\theta\xi^\alpha_1})| \ge g(1)/2.$$

Consider

$$g(\delta) := \inf_{\alpha\in[0,\alpha^*]} \inf_{|\theta|\ge\delta} |1 - E^\alpha_\pi(e^{i\theta\xi^\alpha_1})|;$$

we want to show that $g(\delta) > 0$ for all $\delta > 0$. If $\delta \ge 1$, then $g(\delta) \ge g(1)/2 > 0$. Suppose $0 < \delta < 1$. To show $g(\delta) > 0$, it suffices to show

$$\inf_{\alpha\in[0,\alpha^*]} \inf_{\delta\le|\theta|\le 1} |1 - E\alpha_\pi(e^{i\theta\xi^\alpha_1})| > 0.$$

However, since $E^\alpha_\pi(e^{i\theta\xi^\alpha_1}) = E^\alpha_\pi(e^{(\alpha+i\theta)\xi_1})/E^\alpha_\pi(e^{\alpha\xi_1})$, it is easy to see that $E^\alpha_\pi(e^{i\theta\xi^\alpha_1})$ is a continuous function of $(\alpha,\theta)$ for $\alpha \in \Gamma$ and real $\theta$. Therefore, $|1 - E^\alpha_\pi(e^{i\theta\xi^\alpha_1})|$, being a continuous function on the compact set $\{(\alpha,\theta) : 0 \le \alpha \le \alpha^*, \delta \le |\theta| \le 1\}$, must attain its minimum there. To complete the proof, we need only show that this minimum value cannot be 0. Supposing to the

contrary that $E_\pi^\alpha(e^{i\theta\xi_1^\alpha}) = 1$ for some $0 \le \alpha \le \alpha^*$ and $\delta \le |\theta| \le 1$, we would have

$$P_\pi^\alpha\left\{\frac{\theta}{2\pi}\xi_1 \text{ is an integer }\right\} = 1.$$

However, by the assumption that $\{(X_n, S_n), n \ge 0\}$ is strongly nonlattice and the property of exponential embedding that $P_\pi$ is absolutely continuous with respect to $P_\pi^\alpha$, this is a contradiction. □

Use the same notation as the paragraph before Theorem 2 in Section 2. Note that $\tau_n$, $\tau_+$ and $\tau_-$ depend on $\alpha$; we omit it here for simplicity. The following lemma is related to uniform strong nonlattice of the ladder chains. The proof is a straightforward generalization of Theorem 6 in Fuh and Lai (1998) and is omitted.

LEMMA 13. *Assume the conditions of Theorem 3 hold. Let $P_{\pi_+}^\alpha$ be the transition probability of the ladder Markov chain $\{(X_{\tau_n}^\alpha, S_{\tau_n}^\alpha), n \ge 0\}$. Then, there exists $\alpha^* > 0$ such that, for $0 \le \alpha \le \alpha^*$, the family $\{(X_{\tau_n}^\alpha, S_{\tau_n}^\alpha), n \ge 0\}$ is uniformly strong nonlattice.*

The following lemma generalizes Lemmas 4.4 and 4.5 in Heyde (1964) for simple random walks.

LEMMA 14. *Assume the conditions of Theorem 3 hold. Then, there exist $\alpha^* > 0$, $r_1 > 0$ and $C$ such that, for all $\alpha \in [0, \alpha^*]$,*

$$E_\pi^\alpha(e^{r_1 S_{\tau_+}^\alpha}) \le C.$$

PROOF. Under the assumption C6 and Lemma 11, for $\alpha \in \Gamma \subset \mathbf{R}$, we can define the linear operators $\mathbf{P}_\alpha$, $\mathbf{P}$, $\nu_*$ and $\mathbf{Q}$ on $\mathcal{N}$ as in (2.11). By the spectral decomposition theory for linear operator on the space $\mathcal{N}$ developed in Proposition 1, we have for $h \in \mathcal{N}$,

$$E_\nu\{e^{\alpha S_n} h(X_n)\} = \lambda^n(\alpha)\nu_*\mathbf{Q}_\alpha h + \nu_*\mathbf{P}_\alpha^n(I - \mathbf{Q}_\alpha)h, \tag{6.1}$$

where $\mathbf{Q}_\alpha$ is defined in (2.12). It also can be shown that there exist $K^* > 0$ and $0 < \delta^* < \delta$ such that, for $|\alpha| \le \delta^*$,

$$\|\nu_*\mathbf{P}_\alpha^n(I - \mathbf{Q}_\alpha)h\| \le K^*\|h\||\alpha|\{(1 + 2\rho)/3\}^n, \tag{6.2}$$

and under assumption C6, it follows from Proposition 1 that $\lambda(\alpha)$ has the Taylor expansion

$$\lambda(\alpha) = 1 + \sum_{j=1}^{r} \lambda_j \alpha^j / j! + \Delta(\alpha) \tag{6.3}$$

in some neighborhood of the origin, where $\Delta(\alpha) = O(|\alpha|^r)$ as $\alpha \to 0$.

Now, for such $\alpha$, we have for all $-\infty < s < \infty$,

$$P_\nu(S_n \le s) \le e^{-\alpha s}(\lambda^n(\alpha)\nu_* \mathbf{Q}_\alpha 1 + \nu_* \mathbf{P}_\alpha^n (I - \mathbf{Q}_\alpha)1),$$

where 1 denotes the identity function. Also, under the assumptions $\mu < 0$ and C6, there exists sufficiently small $\alpha$ such that $\lambda(\alpha) < 1$. Along this with (6.2), there exists $C > 0$ such that $e^C \lambda(\alpha) < 1$ and for all $c$, $0 < c < C$,

$$\sum_{n=1}^{\infty} e^{cn} P_\nu(S_n \le s) < \infty. \tag{6.4}$$

Next, for all $\gamma \in (0,1)$, define $F_0(s) = I_{\{s \ge 0\}}$, $F_1(s) = P_\nu(S_1 \le s)$ and $F_n(s) = P_\nu(S_n \le s; \max_{1 \le k \le n-1} S_k \le \log \gamma)$, for $n > 1$. Then, (6.4) implies that

$$\sum_{n=1}^{\infty} e^{rn} F_n(\log \gamma) < \infty, \tag{6.5}$$

for some $r > 0$.

Note that the probability $p_n$ of the first passage time $\tau(\gamma)$ out of the interval $(\log \gamma, \infty)$ for the Markovian random walk $S_n$ is $n$ is given by

$$p_n = F_{n-1}(\log \gamma) - F_n(\log \gamma), \qquad n \ge 1. \tag{6.6}$$

By (6.5) and (6.6), we have $E_\nu e^{t\tau(\gamma)} < \infty$ for some $t > 0$, for all $\gamma \in (0,1)$. Hence

$$E_\pi(e^{t\xi_1}) < \infty \quad \text{implies} \quad E_\pi(e^{tS_{\tau_+}}) < \infty. \tag{6.7}$$

Using the requirement in the definition of the exponential embedding that $\Gamma$ must contain an interval about 0, take any positive $\alpha_1 \in \Gamma$. Let $C := E_\pi(e^{\alpha_1 S_{\tau_+}})$; by (6.7), $C$ is finite. Since $0 < \inf_{|\alpha| > \delta, x \in \mathcal{X}} r(x;\alpha) \le \sup_{|\alpha| > \delta, x \in \mathcal{X}} r(x;\alpha) < \infty$, then, if we take $\alpha^*$ and $r_1$ both to be $\alpha_1/2$, say, for any $\alpha \in [0, \alpha^*]$ we have

$$\begin{aligned} E_\pi^\alpha\{e^{r_1 S_{\tau_+}^\alpha}\} &= E_\pi\left\{\frac{r(X_{\tau_+};\alpha)}{r(X_0;\alpha)} e^{(r_1+\alpha)S_{\tau_+} - \tau_+ \psi(\alpha)}\right\} \\ &\le E_\pi\left\{\frac{r(X_{\tau_+};\alpha)}{r(X_0;\alpha)} e^{(r_1+\alpha)S_{\tau_+}}\right\} \\ &\le E_\pi\{e^{\alpha_1 S_{\tau_+}}\} = C, \end{aligned}$$

which is the desired property. $\square$

LEMMA 15. *Assume the conditions of Theorem 3 hold. Suppose $b \to \infty$ and $0 < \alpha \downarrow 0$ such that, for some $-\infty < \vartheta < \infty$, $\alpha b \to \vartheta$. Then for $0 \le t, s \le \infty$:*

(i) $P_\pi^\alpha\{\tau(b) \le b^2 t, S^\alpha_{\tau(b)} - b \le s\} \longrightarrow G(t;\vartheta,1)H_+(s)$ *and*

(ii) $E_\pi^\alpha l\{(S^\alpha_{\tau(b)} - b)^r; \tau(b) < \infty\} \longrightarrow \frac{E_{\pi_+} S^{r+1}_{\tau_+}}{(r+1)E_{\pi_+} S_{\tau_+}}$.

PROOF. (i) Let $m = b^2 t$ and $\mathcal{F}_N$ be the $\sigma$-algebra generated by $\{(X_n, S_n), n \le N\}$. By Wald's likelihood ratio identity for Markov chains, we have that for any stopping time $N$, $\alpha', \alpha'', A \in \mathcal{F}_N$, and for each fixed $x \in \mathcal{X}$,

$$
\begin{aligned}
&P_x^{\alpha'}\{A \cap (N < \infty)\} \\
&\quad = \int_{A \cap (N<\infty)} \frac{r(X_N;\alpha')}{r(x;\alpha')} \exp((\alpha' - \alpha'')S_N - N(\Lambda(\alpha') - \Lambda(\alpha''))) \, dP_x^{\alpha''}.
\end{aligned} \tag{6.8}
$$

And this implies that

$$
\begin{aligned}
&P_\pi^\alpha\{\tau \le m, S_\tau^\alpha - b \le s\} \\
&\quad = E_\pi\left[\frac{r(X_\tau;\alpha)}{r(X_0;\alpha)} \exp\{\alpha S_\tau - \tau\Lambda(\alpha)\}; \tau \le m, S_\tau - b \le s\right] \\
&\quad = \exp(\alpha b) E_\pi\left[\frac{r(X_\tau;\alpha)}{r(X_0;\alpha)} \exp\{\alpha(S_\tau - b) - \tau\Lambda(\alpha)\};\right. \\
&\qquad\qquad\qquad\qquad \left. \tau \le m, S_\tau - b \le s\right].
\end{aligned} \tag{6.9}
$$

It follows that as $0 < \alpha \downarrow 0$, $\alpha b \to \vartheta$, $r(X_\tau;\alpha)/r(X_0;\alpha) \to 1$ and $\Lambda(\alpha) \sim \frac{1}{2}\alpha^2 \sim \frac{1}{2}\vartheta^2/b^2$. Hence at least for all finite $s$, Lemma 7 shows that the right-hand side of (6.8) converges to

$$
\begin{aligned}
&\exp(\vartheta) E_\pi[\exp\{-\tfrac{1}{2}\vartheta^2 \tau_W(1)\}; \tau_W(1) \le t] H_+(s) \\
&\quad = E_\pi[\exp\{\vartheta W(\tau_W(1)) - \tfrac{1}{2}\vartheta^2 \tau_W(1)\}; \tau_W(1) \le t] H_+(s) \\
&\quad = P_\pi^{(\vartheta)}\{\tau_W(1) \le t\} H_+(s) \\
&\quad = G(t;\vartheta,1) H_+(s).
\end{aligned}
$$

That this calculation is also valid when $s = \infty$ follows from the Markov renewal Theorem 1 once it is known that $E_\pi(e^{rS_{\tau_+}}) < \infty$, for some $r > 0$. This holds by Lemma 14.

(ii) The proof of the convergence of $E_{\pi_+} S_{\pi_+}$ follows from Lemma 6. The rest is a similar calculation and is omitted. □

By using the exponential martingale (2.13), uniform renewal theorem in Theorem 1 and Lemmas 11–15, the proof of Theorem 3 is similar to that of Theorem 2 and is omitted.

# APPENDIX

**Characteristic functions of uniform Markov random walks.** Here we generalize the work of Fuh and Lai (2001). We include it for completeness. Using the same notation and assumptions as in the first paragraph after Theorem 2 of Section 2, define $\mathbf{P}_z^\alpha$, $\mathbf{P}$, $\nu_*^\alpha$ and $\mathbf{Q}$ on $\mathcal{N}$ as (2.11). Condition K2 ensures that $\mathbf{P}_z^\alpha$ and $\mathbf{P}$ are bounded linear operators on $\mathcal{N}$, and (2.2) implies that

$$\|\mathbf{P}^n - \mathbf{Q}\| = \sup_{h \in \mathcal{N}:\, \|h\|=1} \|\mathbf{P}^n h - \mathbf{Q}h\| \le \gamma \rho^n. \tag{A.1}$$

For a bounded linear operator $\mathbf{T}:\mathcal{N} \to \mathcal{N}$, the resolvent set is defined as $\{y \in \mathbf{C}: (\mathbf{T} - yI)^{-1} \text{ exists}\}$ and $(\mathbf{T} - yI)^{-1}$ is called the resolvent (when the inverse exists). From (A.1) it follows that, for $y \ne 1$ and $|y| > \rho$,

$$R(y) := \mathbf{Q}/(y-1) + \sum_{n=0}^{\infty} (\mathbf{P}^n - \mathbf{Q})/y^{n+1} \tag{A.2}$$

is well defined. Since $R(y)(\mathbf{P} - yI) = -I = (\mathbf{P} - yI)R(y)$, the resolvent of $\mathbf{P}$ is $-R(y)$. Moreover, by K3 and an argument similar to the proof of Lemma 2.2 of Jensen (1987), there exist $K > 0$ and $\eta > 0$ such that, for $|z| \le \eta$, $|y-1| > (1-\rho)/6$ and $|y| > \rho + (1-\rho)/6$, $\|\mathbf{P}_z^\alpha - \mathbf{P}\| \le K|\alpha|$ and $R_z^\alpha(y) := \sum_{n=0}^\infty R(y)\{(\mathbf{P}_z^\alpha - \mathbf{P})R(y)\}^n$ is well defined. Since $R_z^\alpha(y)(\mathbf{P}_z^\alpha - yI) = R_z^\alpha(y)\{(\mathbf{P}_z^\alpha - \mathbf{P}) + (\mathbf{P} - yI)\} = -I = (\mathbf{P}_z^\alpha - yI)R_z^\alpha(y)$, the resolvent of $\mathbf{P}_z^\alpha$ is $-R_z^\alpha(y)$.

For $|z| \le \eta$, the spectrum (which is the complement of the resolvent set) of $\mathbf{P}_z^\alpha$ therefore lies inside the two circles $C_1 = \{y: |y-1| = (1-\rho)/3\}$ and $C_2 = \{y: |y| = \rho + (1-\rho)/3\}$. Hence by the spectral decomposition theorem [cf. Riesz and Sz-Nagy (1955), page 421], $\mathcal{N} = \mathcal{N}_1(z) \oplus \mathcal{N}_2(z)$ and

$$\mathbf{Q}_z^\alpha := \frac{1}{2\pi i}\int_{C_1} R_z^\alpha(y)\,dy, \qquad I - \mathbf{Q}_z^\alpha := \frac{1}{2\pi i}\int_{C_2} R_z^\alpha(y)\,dy \tag{A.3}$$

are parallel projections of $\mathcal{N}$ onto the subspaces $\mathcal{N}_1(z)$, $\mathcal{N}_2(z)$, respectively. Moreover, by an argument similar to the proof of Lemma 2.3 of Jensen (1987), there exists $0 < \delta \le \eta$ such that $\mathbf{B}_1(z)$ is one-dimensional for $|z| \le \delta$ and $\sup_{|z|\le\delta} \|\mathbf{Q}_z^\alpha - \mathbf{Q}\| < 1$. For $|z| \le \delta$, let $\lambda^\alpha(z)$ be the eigenvalue of $\mathbf{P}_z^\alpha$ with corresponding eigenspace $\mathcal{N}_1(z)$. Since $\mathbf{Q}_z^\alpha$ is the parallel projection onto the subspace $\mathbf{B}_1(z)$ in the direction of $\mathbf{B}_2(z)$, (2.12) holds. Therefore, for $h \in \mathcal{N}$,

$$\begin{aligned} E_\nu\{e^{zS_n^\alpha} h(X_n)\} &= \nu_*^\alpha \mathbf{P}_z^{\alpha\, n} h = \nu_*^\alpha \mathbf{P}_z^{\alpha\, n}\{\mathbf{Q}_z^\alpha + (I - \mathbf{Q}_z^\alpha)\}h \\ &= (\lambda^\alpha(z))^n \nu_*^\alpha \mathbf{Q}_z^\alpha h + \nu_*^\alpha \mathbf{P}_z^{\alpha\, n}(I - \mathbf{Q}_z^\alpha)h. \end{aligned}$$

Suppose K4 also holds. An argument similar to the proof of Lemma 2.4 of Jensen (1987) shows that, there exist $0 < \delta^* < \delta$ and $K^* > 0$ such that for $|z| \le \delta^*$, $|\nu_*^\alpha \mathbf{P}_z^{\alpha\, n}(I - \mathbf{Q}_z^\alpha)h| \le K^* \|h\|_w |z| \{(1+2\rho)/3\}^n$. Moreover, analogous to Lemmas 2.5–2.7 of Jensen (1987), it can be shown that $\lambda^\alpha(z)$, $\nu_*^\alpha \mathbf{Q}_z^\alpha h$ and $\sum_{n=0}^\infty \nu_*^\alpha \mathbf{P}_z^{\alpha\, n}(I - \mathbf{Q}_z^\alpha)h$ have continuous partial derivatives of order $[r]$ for $|z| \le \delta^*$. Furthermore, we have the following proposition.

Proposition 1. *Assume the conditions of Theorem* 1 *hold. Let* $h \in \mathcal{N}$ *and there exists a* $\delta > 0$ *such that* $z \in \mathbf{C}$ *and* $|z| \le \delta$.

(i) $E_\nu\{e^{zS_n^\alpha} h(X_n)\} = (\lambda^\alpha(z))^n \nu_*^\alpha \mathbf{Q}_z^\alpha h + \nu_*^\alpha (\mathbf{P}_z^\alpha)^n (I - \mathbf{Q}_z^\alpha)h$. *Moreover, there exist* $0 < \delta^* < \delta$, $0 < \gamma < 1$ *and* $K > 0$ *such that, for* $|z| \le \delta^*$, $\lambda^\alpha(z)$, $\nu_*^\alpha \mathbf{Q}_z^\alpha h$ *and* $\sum_{n=0}^\infty \nu_*^\alpha (\mathbf{P}_z^\alpha)^n (I - \mathbf{Q}_z^\alpha)h$ *have continuous partial derivatives of order* $[r]$, *and*

$$|\nu_*^\alpha (\mathbf{P}_z^\alpha)^n (I - \mathbf{Q}_z^\alpha)h| \le K\|h\||z|\gamma^n \qquad \text{for all } n \ge 1.$$

*Furthermore,*

$$\lambda^\alpha(0) = 1, \qquad \nabla\lambda^\alpha(0) = i\Gamma\mu^\alpha, \qquad \nabla^2\lambda^\alpha(0) = -\Gamma V^\alpha \Gamma'.$$

(ii) *Define* $f_A^\alpha(z) = \sum_{n=0}^\infty E_\nu(e^{zS_n^\alpha} \mathbb{1}_{\{X_n \in A\}})$, *and let* $h_A(z) = \mathbb{1}_{\{x \in A\}}$. *Then for* $0 < |z| \le \delta^*$,

$$f_A^\alpha(z) = (1 - \lambda^\alpha(z))^{-1} \nu_*^\alpha \mathbf{Q}_z^\alpha h_A + \eta^\alpha(z),$$

*where* $\eta^\alpha(z)$ *has continuous partial derivatives of order* $[r]$ *and* $\eta^\alpha(z) = O(|z|)$ *as* $z \to 0$.

**Acknowledgments.** The author is grateful to the referee for constructive comments, suggestions and correction of some mistakes in the proofs. Thanks also to Professor Søren Asmussen, the Associate Editor, for careful proofreading.

Institute of Statistical Science
Academia Sinica
Taipei 11529, Taiwan
Republic of China
e-mail: stcheng@stat.sinica.edu.tw